 \documentclass[3p,sort&compress]{elsarticle}
\usepackage[cp1252]{inputenc}
\usepackage{textcomp}
\usepackage{amsmath}
\usepackage{amssymb}
\usepackage{graphicx}
\usepackage[abs]{overpic}
\usepackage{hyphenat}
\usepackage[ruled,vlined]{algorithm2e}
\usepackage{inputenc}
\usepackage{amsthm}
\usepackage{booktabs}
\usepackage{bbold}
\usepackage{enumitem}
\setenumerate[1]{label = \arabic*.}
\setenumerate[2]{label*= \arabic*}
\setenumerate[3]{label*=.\arabic*}
\setenumerate[4]{label*=.\arabic*}
\setenumerate[4]{label*=.\arabic*}
\usepackage{xcolor}
\usepackage{multirow}
\usepackage[titletoc]{appendix}
\usepackage{anyfontsize}
\usepackage{framed}
\usepackage{footnote}
\usepackage{mathbbol}
\usepackage{subcaption}
\usepackage{tabulary}
\usepackage{bbm} 
\newtheoremstyle{dotless}{}{}{\itshape}{}{\bfseries}{}{ }{}
\theoremstyle{dotless}

\makeatletter

\renewcommand{\mathbf}{\boldsymbol}

\graphicspath{{figures/}}

\usepackage{color}
\usepackage{soul}
\usepackage{amsthm}

\usepackage{units}

\journal{\vspace*{1mm}}

\makeatother

\newcommand{\be}{\begin{equation}}
\newcommand{\ee}{\end{equation}}

\newcommand{\bC}{\mbox{\boldmath{$C$}}}
\newcommand{\bd}{\mbox{\boldmath{$d$}}}

\newcommand{\mbe}{\mbox{\boldmath{$e$}}}
\newcommand{\fb}{\mbox{\boldmath{$f$}}}
\newcommand{\Fb}{\mbox{\boldmath{$F$}}}

\newcommand{\bK}{\mbox{\boldmath{$K$}}}

\newcommand{\bn}{\mbox{\boldmath{$n$}}}
\newcommand{\bN}{\mbox{\boldmath{$N$}}}

\newcommand{\bP}{\mbox{\boldmath{$P$}}}

\newcommand{\bt}{\mbox{\boldmath{$t$}}}

\newcommand{\bR}{\mbox{\boldmath{$R$}}}
\newcommand{\bu}{\mbox{\boldmath{$u$}}}
\newcommand{\bU}{\mbox{\boldmath{$U$}}}

\newcommand{\bv}{\mbox{\boldmath{$v$}}}

\newcommand{\bx}{\mbox{\boldmath{$x$}}}
\newcommand{\bX}{\mbox{\boldmath{$X$}}}

\newcommand{\bphi}{\mbox{\boldmath{$\phi$}}}
\newcommand{\bvarphi}{\mbox{\boldmath{$\varphi$}}}

\newcommand{\bxi}{\mbox{\boldmath{$\xi$}}}

\newcommand{\TEfive}{\mbox{$\text{T}_{5}(\lambda)$}}
\newcommand{\TEn}{\mbox{$\text{T}_{n}(\lambda)$}}

\newfont{\twelvemsb}{msbm10 at 11.6pt}

\newcommand{\Div}{{\rm Div}}

\newcommand{\Grad}{{\rm Grad}}

\newcommand{\tr}{\mathop{\rm tr}}

\begin{document}
\begin{frontmatter}

\title{A virtual element method for isotropic hyperelasticity}

\author[danieladdress]{D.~van~Huyssteen\corref{cor2}}
\author[reddyaddress]{B.~D.~Reddy}

\address[danieladdress]{Department of Mechanical Engineering and Centre for Research in Computational and Applied Mechanics, University of Cape Town, Rondebosch, 7701, South Africa}

\address[reddyaddress]{Department of Mathematics and Applied Mathematics and Centre for Research in Computational and Applied Mechanics, University of Cape Town, Rondebosch, 7701, South Africa}



 \cortext[cor2]{Corresponding author. Tel.:
 +27 (0) 21-650-3817; E-mail address:
 vhydan001@myuct.ac.za}%

\begin{abstract}
This work considers the application of the virtual element method to plane hyperelasticity problems with a novel approach to the selection of stabilization parameters. The method is applied to a range of numerical examples and well known strain energy functions, including neo-Hookean, Mooney-Rivlin and Ogden material models. For each of the strain energy functions the performance of the method under varying degrees of compressibility, including near-incompressibility, is investigated. Through these examples the convergence behaviour of the virtual element method is demonstrated. Furthermore, the method is found to be robust and locking free for a variety of element geometries, including elements with a high degree of concavity.
\end{abstract}


\end{frontmatter}

\section{Introduction}\label{S:Introduction}
The virtual element method (VEM) is a recent extension of the well established finite element method for approximating solutions to problems posed as systems of partial differential equations, or in their variational form \cite{VEIGA2012,Veiga2014}. 

Several adaptations of the finite element method have been developed to overcome specific challenges. Mixed methods allow all variables of interest to be approximated explicitly and have been used successfully for problems involving near-incompressibility, and those in which the geometry is characterised by a small length scale, characteristics that, for low-order finite elements, lead to volumetric and shear locking respectively \cite{Boffi2013,Hughes2000}. Another development is the discontinuous Galerkin (DG) method in which interelement continuity is abandoned, allowing for greater flexibility with regard to meshing \cite{Arnold2002}. Additionally, when designed appropriately, the DG method is stable and uniformly convergent in the case of near-incompressibility for low-order approximations \cite{Grieshaber2015,Hansbo2002,Wihler2004}.

In contrast to the finite element method, in which elements are typically triangular and quadrilateral in 2D or tetrahedral and hexahedral in 3D, the virtual element method permits arbitrary polygonal or polyhedral elements in 2D and 3D respectively. The inherent flexibility of virtual element meshes lends the method to problems involving complex geometries and adaptive meshing. 

The virtual element method has been applied to a growing range of problems in solid mechanics, including isotropic \cite{Artioli2017} and anisotropic \cite{Reddy2019} linear elasticity as well as their non-linear counterparts \cite{Chi2017,WriggersIsotropic2017,WriggersAnisotropic2017}. The method has also been applied to problems involving viscoelasticity and shape memory alloys in \cite{Artioli2017a}, and to plasticity \cite{Artioli2017a,WriggersPlastic2017,HudobivnikPlastic2018}, coupled thermo-elasticity \cite{Dhanush2018} and thermo-plasticity \cite{Aldakheel2019}. Further applications include scalar damage models \cite{Bellis2018} as well as brittle \cite{Aldakheel2018} and ductile \cite{Aldakheel2019a} fracture models.

To date the bulk of the literature focuses on low-order virtual elements; however, there is increasing interest in higher order serendipity-like virtual element methods \cite{Veiga2016,Veiga2016a,Veiga2018,Bellis2019}.

The virtual element method may be characterised by the splitting of the field describing the variable of interest $\bphi$ into a consistent part, computed via a projection operator $\Pi (\bphi)$, and an error-like part defined by the difference between the field and its projection $\bphi - \Pi (\bphi)$. Use of only the consistent term would result in a rank deficient stiffness matrix, thus necessitating the introduction of a stabilization term. In the linear case the stabilization term can be easily computed from a sum \cite{Gain2014,Veiga2015,Veiga2013} or product \cite{Artioli2017,Reddy2019} of nodal values multiplied by a scalar stabilization parameter representative of the material's constitution. Additionally, it was shown computationally in \cite{Reddy2019} for the case of transversely isotropic materials that, the stabilization parameter can be chosen in such a way that a low-order virtual element method is locking free in the case of near-incompressibility.

In the non-linear case it is also possible to base the stabilization term on a sum of nodal values. This method was successfully implemented in \cite{Chi2017}, with the stabilization parameter computed from the fourth-order elasticity tensor, for problems with small load steps. The small load steps are necessary as the stabilization parameter is computed from the deformation of the previous load step. This method is, however, inefficient as a result of the constraint on load step size. 

A popular approach to computing the stabilization term involves defining a stabilization strain energy function $\hat{\Psi}$ \cite{WriggersIsotropic2017} and approximating the difference $\Psi(\bphi)-\Psi(\Pi(\bphi))$ by $\hat{\Psi}(\bphi)-\hat{\Psi}(\Pi(\bphi))$, in which the integration of $\hat{\Psi}(\bphi)$ is performed using a quadrature rule over sub-elements. It is clear then that if the material strain energy function $\Psi(\Pi(\bphi))$ is the same as the stabilization strain energy function the method would reduce to an inscribed mesh of the sub-elements. To this end, modified Lam\'{e} parameters are typically used for the stabilization strain energy function with these parameters typically computed in one of two ways. The first method involves the assumption of a compressible Poisson's ratio, which prevents locking, and computation of a scalar value that is a function of the geometry of the virtual element. Specifically, the convex hull of the element is required. This value of Poisson's ratio and the geometry factor are then used to compute the modified Lam\'{e} parameters; see for example \cite{WriggersAnisotropic2017,WriggersPlastic2017,HudobivnikPlastic2018}. For the second method the modified Lam\'{e} parameters are computed by multiplying the standard parameters by a chosen scalar value $\gamma$ in the range $(0,1]$, for which $\gamma=1$ results in the degenerate case of inscribed sub-elements \cite{Bellis2019,Aldakheel2018}.

To the best of the authors' knowledge only neo-Hookean models have been used to describe the material material behaviour and the stabilization strain energy in works dealing with non-linear elasticity.

In this work we present an alternative approach to stabilization that utilizes a different geometry factor and the introduction of an incompressibility factor to scale the Lam\'{e} parameter $\mu$. Additionally, the method uses a truncated Taylor expansion of the volumetric Lam\'{e} parameter $\lambda$, and does not require the assumption of a compressible Poisson's ratio. We investigate the effect of compressibility, determined by Poisson's ratio, on the proposed stabilization, the effect of element geometry, including highly concave elements; as well as non star-shaped elements.

We also extend earlier investigations by considering a range of well known strain energy functions.

The structure of the rest of this work is as follows. Section \ref{sec:GovEq} sets out the kinematics and constitutive relations governing hyperelasticity and the associated variational form. The details of the virtual element method and the proposed stabilization methodology are presented in Section \ref{sec:VEM}, and the set of numerical results are presented in Section \ref{sec:Results}. The work concludes with a summary of results.

\section{Governing equations of hyperelasticity} \label{sec:GovEq}
Consider an elastic body occupying a plane, polygonal bounded domain $\Omega \subset \mathbb{R}^{2} $ with boundary $\partial \Omega$ in its undeformed configuration. The boundary comprises a non-trivial Dirichlet part $\Gamma_{D}$ and a Neumann part $\Gamma_{N}$ such that $\Gamma_{D} \cap \Gamma_{N} = \emptyset$ and $\overline{\Gamma_{D} \cup \Gamma_{N}}=\partial \Omega$.
\begin{figure}[ht!]
	\centering
	\def\svgwidth{0.4\textwidth}
	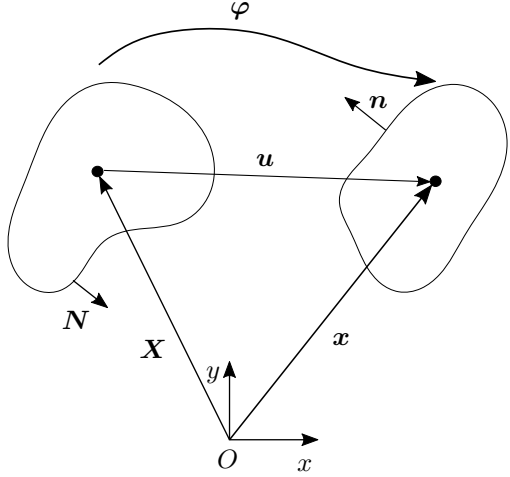
	\vspace*{2mm}
	\caption{Motion of a body}
	\label{fig:Potato}
\end{figure}
\newline \noindent
The body undergoes a motion $\bvarphi$, see Figure \ref{fig:Potato}, such that the current position of a point $\bx$, initially at $\bX$, is given by
\begin{align}
\bx &= \bvarphi(\bX,t) \nonumber \\
&= \bX + \bu(\bX,t) \, ,
\end{align}
where $\bu$ is the displacement. We denote by $\Fb$ the deformation gradient defined by
\begin{align}
\Fb &= \Grad \, \bvarphi \nonumber \\
&= \boldsymbol{1} + \Grad \, \bu \, ,
\end{align}
where $\Grad (\bullet){ij}=\nabla(\bullet)=\frac{\partial}{\partial X_{j}}(\bullet)_{i}$. For equilibrium we require
\begin{equation}
- \Div \, \bP = \fb \, ,
\end{equation}
where $\fb$ denotes the body force, $\bP$ the first Piola-Kirchhoff stress and $\Div (\bullet)=\nabla\cdot(\bullet)=\frac{\partial}{\partial X_{i}}(\bullet)\cdot\mbe_{i}$. The Dirichlet and Neumann boundary conditions are given by
\begin{align}
&\bu = \bar{\bu} \quad \text{on } \Gamma_{D} \, , \\
&\bP \, \bN = \bar{\bt} \quad \text{on } \Gamma_{N} \, ,
\end{align}
where $\bar{\bu}$ denotes a prescribed displacement, $\bar{\bt}$ a surface traction, and $\bN$ the outward unit normal vector. We denote by $\Psi(\bvarphi)$ a strain energy function from which the first Piola-Kirchhoff stress can be found from
\begin{equation}
\bP=\frac{\partial \Psi}{\partial \Fb} \, .
\end{equation}
The strain energy functions considered in this work are functions of the invariants of the right Cauchy-Green tensor defined by
\begin{equation}
\bC  = \Fb^{T} \Fb \, ,
\end{equation}
and the principal stretches $\lambda_{i}$, which are the eigenvalues of the stretch tensor $\bU=\bC^{1/2}$.
The invariants of $\bC$ are given by
\begin{align}
I_{C}=\tr \bC \quad &\text{or} \quad I_{C}=\lambda_{1}^{2} + \lambda_{2}^{2} + \lambda_{3}^{2} \, , \\
II_{C}=\frac{1}{2}\left[ \left( \tr \bC \right)^{2} - \tr \bC ^{2} \right] \quad &\text{or} \quad II_{C}= \lambda_{1}^{2} \lambda_{2}^{2} +  \lambda_{2}^{2} \lambda_{3}^{2} + \lambda_{3}^{2} \lambda_{1}^{2} \, \,  \text{and} \\
III_{C}=\det \bC \quad &\text{or} \quad III_{C}=\lambda_{1}^{2} \lambda_{2}^{2} \lambda_{3}^{2} \, .
\end{align}
We note that the case $\lambda_{3}=1$ corresponds to the condition of plane strain. Additionally, the invariants of the isochoric part of the right Cauchy-Green tensor $\overline{\bC}=\left(\det \bC\right)^{-1/3} \bC$ are required, and are given by
\begin{align}
\overline{I}_{C}=\tr \overline{\bC} \quad &\text{or} \quad \overline{I}_{C}=\bar{\lambda}_{1}^{2} + \bar{\lambda}_{2}^{2} + \bar{\lambda}_{3}^{2} \, , \label{eqn:iso1} \\
\overline{II}_{C}=\frac{1}{2}\left[ \left( \tr \overline{\bC} \right)^{2} - \tr \overline{\bC} ^{2} \right] \quad &\text{or} \quad \overline{II}_{C}= \bar{\lambda}_{1}^{2} \bar{\lambda}_{2}^{2} +  \bar{\lambda}_{2}^{2} \bar{\lambda}_{3}^{2} + \bar{\lambda}_{3}^{2} \bar{\lambda}_{1}^{2} \, \,  \text{and} \\
\overline{III}_{C}=\det \overline{\bC} \quad &\text{or} \quad \overline{III}_{C}=\bar{\lambda}_{1}^{2} \bar{\lambda}_{2}^{2} \bar{\lambda}_{3}^{2} \label{eqn:iso3} \, ,
\end{align}
where $\bar{\lambda}_{i}$ denotes the isochoric part of the $i$\textsuperscript{th} principal stretch. Finally, we denote by $J$ the Jacobian of the deformation gradient, so that
\begin{align}
J&=\det \, \Fb \, .
\end{align}
The variational boundary value problem is then one of solving for the displacement by minimizing the total potential energy
\begin{equation}
U(\bu,\, \nabla \bu)=\int_{\Omega}\left[ \Psi (\nabla \bu) - \fb \cdot \bu \right] \,d\Omega - \int_{\Gamma_{N}}\bt \cdot \bu \, d\Gamma \, .
\end{equation}

\section{The virtual element method}
\label{sec:VEM}
\begin{figure}[ht!]
	\centering
	\def\svgwidth{0.4\textwidth}
	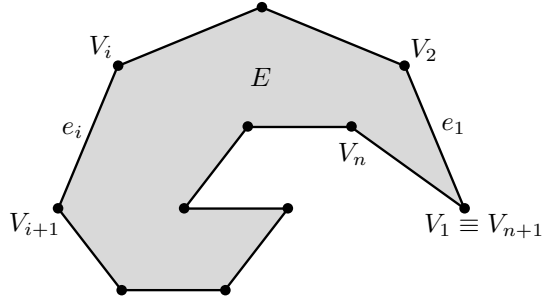
	\caption{Example polygonal element}
	\label{fig:SampleElement}
\end{figure}

The domain $\Omega$ is partitioned into a mesh of elements comprising non-overlapping arbitrary polygons $E$ with $\overline{\cup E}=\bar{\Omega}$. A sample polygonal element is shown in Figure \ref{fig:SampleElement}. We denote by $e_{i}$ the edge connecting vertices $V_{i}$ and $V_{i+1}$ with $i=1,\dots,n_{\rm v}$, where $n_{\rm v}$ is the total number of vertices of element $E$.
 
We construct a conforming approximation in a space $V^{h} \subset V$. The space $V^{h}$ comprises functions that are continuous on $\Omega$, piecewise linear on $\partial \Omega$, and with $\Delta \bv$ vanishing on $E$:
\begin{align}
V^{h}&=\{\bv_{h}\in V \, | \, \bv_{h} \in [C(\Omega)]^{2},\,  \Delta \bv = \boldsymbol{0} \text{ on } E \, ,\, \bv_{h}|_{e} \, \in \, P_{1}(e)  \} \, ,
\end{align}
where $\Delta \bv$ denotes the Laplacian operator, with $ \Delta \bv = \nabla\cdot\nabla\bv$. Here and henceforth $P_{k}(X)$ denotes the space of polynomials of degree $\leq \, k$ on the set $X \, \subset \, \mathbb{R}^{d} $ $(d=1,2)$. We assign degrees of freedom to the nodes, which are located at the element vertices, and write, for each element,
\begin{align}
\bv_{h}|_{E}=\bxi \bd
\end{align}
in which $\bxi$ denotes a matrix of virtual basis functions, and $\bd$ is the $2n_{\rm v} \times 1$ vector of degrees of freedom. All computations are carried out on the edges $e$ of elements, and it is convenient to write
\begin{align}
\bv_{h}|_{\partial E}= \boldsymbol{\mathcal{N}} \bd  , \label{eqn:approximations}
\end{align}
where $\boldsymbol{\mathcal{N}} $ denotes the matrix of standard Lagrangian linear basis functions. Thus the basis functions $\bxi$ are not explicitly known, and not required to be known; their traces on the boundary are however required and are simple Lagrangian functions.
\newline \noindent
We will require the projection $\Pi \, : \, V_{h}|_{E} \, \rightarrow \, P_{0}(E)$, defined on $E$ by
\begin{equation}
|E|\Pi\bv_{h}= \int_{E}\Pi \, \bv_{h} \, dx = \int_{E} \Grad \, \bv_{h} \, dx \, . \label{eqn:AverageProjection}
\end{equation}
Thus $\Pi$ is the $L^{2}$-orthogonal projection onto constants of the gradient of $\bv_{h}$ on an element $E$. From (\ref{eqn:approximations}), and given that $\Pi \, \bv_{h}$ is constant we have, in component form,
\begin{align}
(\Pi \, \bv_{h})_{ij}&=\frac{1}{|E|}\int_{E}(v_{h})_{i,j} \, dx \nonumber \\
&=\frac{1}{|E|}\oint_{\partial E}(v_{h})_{i}N_{j} \, ds \nonumber \\
&=\frac{1}{|E|} \sum_{e\in \partial E}\int_{e} \mathcal{N}_{iA}d_{A}^{E}N_{j}  \, ds \, . \label{eqn:Projection}
\end{align}
Here $d_{A}^{E}$ denotes the degrees of freedom associated with element $E$, summation is implied over all repeated indices, and we have used integration by parts and the representation (\ref{eqn:approximations}). The integrals in (\ref{eqn:Projection}) are readily evaluated as the edge basis functions are known. Thus the projection $\Pi \, \bv_{h}$ is available as a function of the degrees of freedom. 

\begin{figure}[ht!]
	\centering
	\def\svgwidth{0.4\textwidth}
	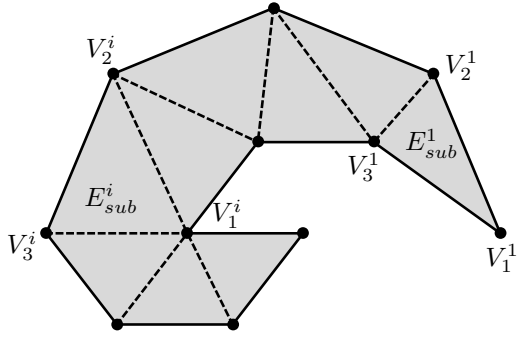
	\caption{Sample decomposed element}
	\label{fig:DecomposedElement}
\end{figure}

To construct the virtual element formulation we start by writing
\begin{align}
U(\bu)&=\sum_{E=1}^{n_{el}}U^{E}(\bu) \, ,
\end{align}
where $n_{el}$ is the total number of elements and $U^{E}(\bu)$ denotes the contribution of element $E$ to the potential energy $U(\bu)$ and is defined by
%
\begin{align}
U^{E}(\bu)&= U_{\Psi}^{E}(\nabla \bu) - U_{f}^{E}(\bu) - U_{t}^{E}(\bu)  \, .
\end{align}
where $U_{\Psi}^{E}(\nabla \bu)$ denotes the internal work and $U_{f}^{E}(\bu)$ and $U_{t}^{E}(\bu)$ denote external work contributions form body $\fb$ and $\bt$ traction forces respectively. The body force contribution is defined by
\begin{equation}
U_{f}^{E}(\bu)=\int_{E} \fb \cdot \bu_{h} \, d\Omega \, ,
\end{equation}
where, as the displacement $\bu_{h}$ on the interior of an element is implicit, the body force term must be approximated in some way. A simple option is to approximate $\fb$ as constant at the element level. The traction contribution is defined by
\begin{equation}
U_{t}^{E}(\bu) = \int_{\Gamma_{N}} \bt \cdot \bu_{h} \, d \Gamma \, ,
\end{equation}
and is computed in the standard way. To compute the internal work contribution we begin by expressing the strain energy in terms of the projection as
\begin{equation}
\Psi(\nabla\bu)=\Psi(\Pi \bu) + \Psi(\nabla \bu) - \Psi(\Pi \bu) \, ,
\end{equation}
which we approximate by
\begin{equation}
\Psi(\nabla\bu)\approx\Psi(\Pi \bu) + \hat{\Psi}(\nabla \bu) - \hat{\Psi}(\Pi \bu) \, ,
\end{equation}
where $\hat{\Psi}$ denotes a modified stabilization strain energy potential defined in the next section. The internal work is then split into a consistency term and a stabilization term
\begin{equation}
U_{\Psi}^{E}(\nabla \bu) = U_{c}^{E} + U_{stab}^{E} \,
\end{equation}
with
\begin{equation}
U_{c}^{E} = \int_{E} \Psi(\Pi \bu_{h}) \, d\Omega \, ,
\end{equation}
which can be computed exactly, and
\begin{equation}
U_{stab}^{E} = \int_{E} \hat{\Psi}(\bu_{h}-\Pi \bu_{h}) \, d\Omega \, , \label{eqn:StabEnergy}
\end{equation}
which must be approximated. As $\bu_{h}$ is only explicitly known on the boundary of an element we choose to perform the integral (\ref{eqn:StabEnergy}) over a set of triangular subdomains \cite{WriggersIsotropic2017}. Figure \ref{fig:DecomposedElement} shows a possible decomposition of an element into $n_{\rm sub}$ triangular subdomains. \\

The computation of (\ref{eqn:StabEnergy}) then becomes
\begin{equation}
U_{stab}^{E}\approx\sum_{i=1}^{n_{\rm sub}}\int_{E_{\rm sub}^{i}}\left[ \hat{\Psi}(\nabla \bu_{h}|_{E_{\rm sub}^{i}}) - \hat{\Psi}(\Pi \bu_{h}|_{E}) \right] \, d \Omega \, ,
\end{equation}
in which the gradient of $\bu_{h}|_{E_{\rm sub}^{i}}$ is assumed to be constant and can be calculated trivially, see \cite{WrigFEM}.
\newline \noindent
All derivations of the potential energy required for the computation of the residual vector $\bR_{E}$ and tangent matrix $\bK_{E}$ for element $E$ are performed using the symbolic tool ACEGEN, see \cite{JozeKorelc2016}. For the consistency term we use the relations 
\begin{align}
\bR_{E}^{c}=\frac{\partial U_{c}(\Pi \bu_{h}|_{E})}{\partial \bd_{E}} \quad \text{and} \quad
\bK_{E}^{c}=\frac{\partial \bR_{E}^{c}(\bd_{E})}{\partial \bd_{E}} \, , 
\end{align}
and for the stabilization term
\begin{align}
\bR_{E}^{stab}=\frac{\partial U_{stab}(\bd_{E})}{\partial \bd_{E}} \quad \text{and} \quad
\bK_{E}^{stab}=\frac{\partial \bR_{E}^{stab}(\bd_{E})}{\partial \bd_{E}} \, ,
\end{align}
with the total residual vector and tangent matrix given by
\begin{align}
\bR_{E}=\bR_{E}^{c}+\bR_{E}^{stab} \quad \text{and} \quad \bK_{E}=\bK_{E}^{c}+\bK_{E}^{stab} \, .
\end{align}

\subsection{Stabilization}
\label{subsec:Stab}
We propose a stabilization strain energy function based on a neo-Hookean form, that is,
\begin{equation}
\hat{\Psi} = \frac{\hat{\mu}}{2}(I_{C} - 3 - 2 \ln (J)) + \frac{\hat{\lambda}}{2}(J-1)^{2} \, ,
\end{equation}
where $\hat{\mu}$ and $\hat{\lambda}$ denote modified Lam\'{e} parameters. The $\hat{\lambda}$ term is modified such that it closely matches the behaviour of $\lambda$ while remaining bounded on $\nu \in \left(-1 \, , \, 0.5\right)$. To this end we employ a fifth-order Taylor expansion of $\lambda$. It is the boundedness of $\TEfive$ that is exploited to prevent volumetric locking of the VEM in the nearly-incompressible limit $\nu \rightarrow 0.5$.

In this work we choose to express material parameters in terms of the familiar engineering constants Young's modulus $E_{y}$ and Poisson's ratio $\nu$ with
\begin{equation}
\lambda = \frac{E_{y} \nu}{(1+\nu)(1-2\nu)} \quad \text{and} \quad 
\mu = \frac{E_{y}}{2(1+\nu)} \, .
\end{equation}
The $n$-th order Taylor expansion of $\lambda$ about $\nu_{0}$ is given by
\begin{equation}
\TEn = \lambda(\nu_{0}) + \sum_{i=1}^{n}\left. \frac{\partial^{i}\lambda}{\partial \nu^{i}} \right|_{\nu_{0}} \frac{(\nu-\nu_{0})^{i}}{i!} \, .
\end{equation}
For simplicity we choose to expand about $\nu_{0}=0^{*}$, which corresponds to a Maclaurian series.

\textsuperscript{*}The choice of $\nu_{0}=0$ is well suited to problems with $\nu \in (-0.5,\, 0.5)$ and prevents locking of the VEM in the incompressible limit. As $\TEfive$ is approximately an odd function about $\nu_{0}$ the opposing limit $\nu \rightarrow -1$ is treated more severely, with $\left| \TEfive|_{\nu = -1} / \TEfive|_{\nu = 0.5} \right| \approx 5 $, this can cause `over stabilization' which manifests as comparatively poor convergence behaviour for coarse meshes. We therefore expect less accurate approximations for $\nu \leq 0$. This, however, is considered acceptable as the vast majority of engineering applications consider values of $\nu \geq 0$. For the most general case one could rather choose $\nu_{0}=-0.25$. Alternatively, a lower-order Taylor expansion could be used.

\begin{figure}[ht!]
	\centering
	\includegraphics[width=0.35\textwidth]{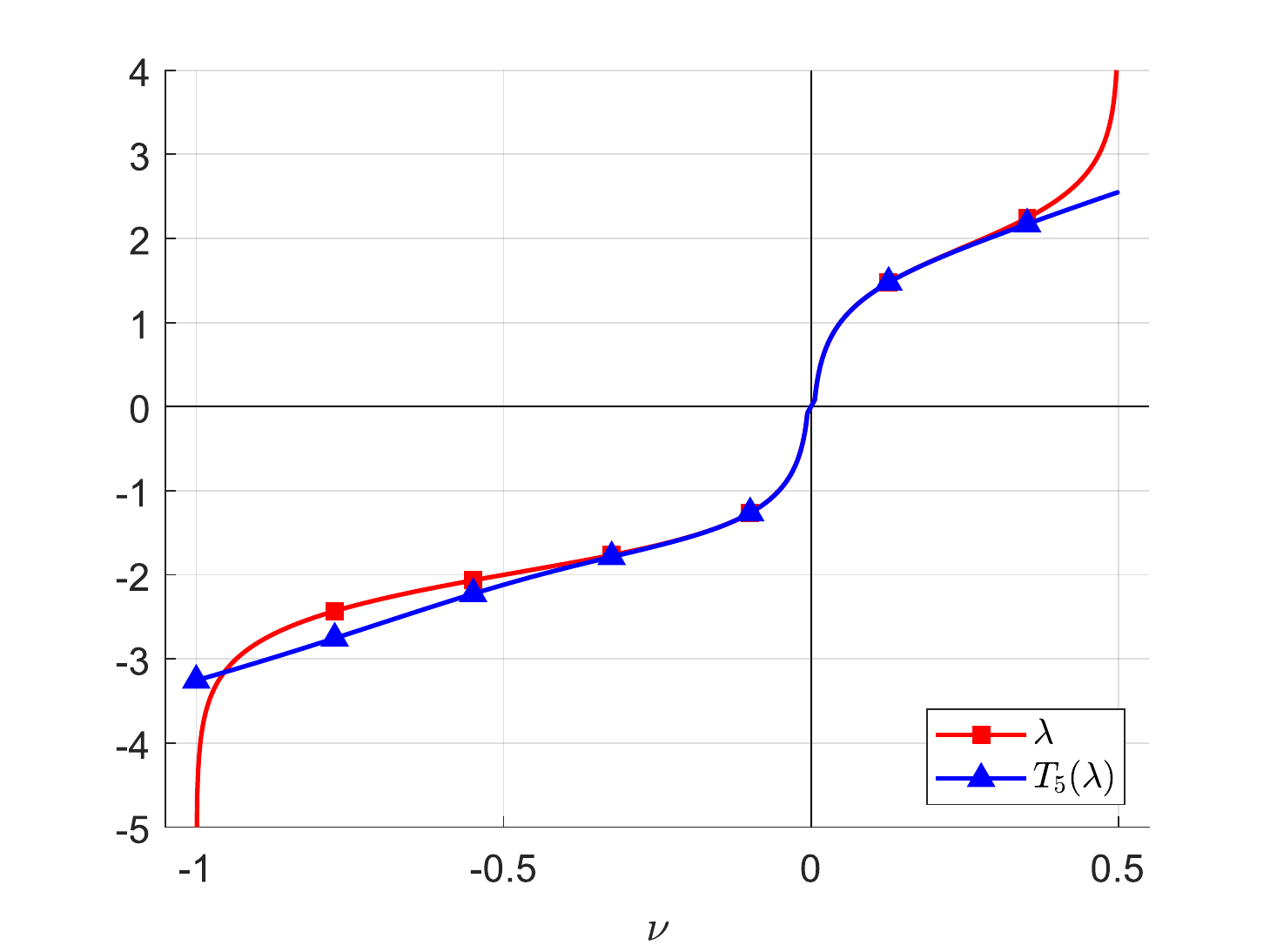}
	\caption{\normalsize Lam\'{e} parameter $\lambda$ and $\TEfive$ vs Poisson's ratio \label{fig:Lame_vs_Poisson}}
\end{figure}

Figure \ref{fig:Lame_vs_Poisson} shows a plot of the Lam\'{e} parameter $\lambda$ and $\TEfive$ vs Poisson's ratio on a logarithmically scaled $y$-axis. We note the boundedness of the Taylor expansion and conversely the divergence of the Lam\'{e} parameter for limiting values of Poisson's ratio.

\begin{figure}[ht!]
	\centering
	\def\svgwidth{0.35\textwidth}
\begingroup%
  \makeatletter%
  \providecommand\color[2][]{%
    \errmessage{(Inkscape) Color is used for the text in Inkscape, but the package 'color.sty' is not loaded}%
    \renewcommand\color[2][]{}%
  }%
  \providecommand\transparent[1]{%
    \errmessage{(Inkscape) Transparency is used (non-zero) for the text in Inkscape, but the package 'transparent.sty' is not loaded}%
    \renewcommand\transparent[1]{}%
  }%
  \providecommand\rotatebox[2]{#2}%
  \newcommand*\fsize{\dimexpr\f@size pt\relax}%
  \newcommand*\lineheight[1]{\fontsize{\fsize}{#1\fsize}\selectfont}%
  \ifx\svgwidth\undefined%
    \setlength{\unitlength}{503.06918899bp}%
    \ifx\svgscale\undefined%
      \relax%
    \else%
      \setlength{\unitlength}{\unitlength * \real{\svgscale}}%
    \fi%
  \else%
    \setlength{\unitlength}{\svgwidth}%
  \fi%
  \global\let\svgwidth\undefined%
  \global\let\svgscale\undefined%
  \makeatother%
  \begin{picture}(1,1)%
    \lineheight{1}%
    \setlength\tabcolsep{0pt}%
    \put(0,0){\includegraphics[width=\unitlength,page=1]{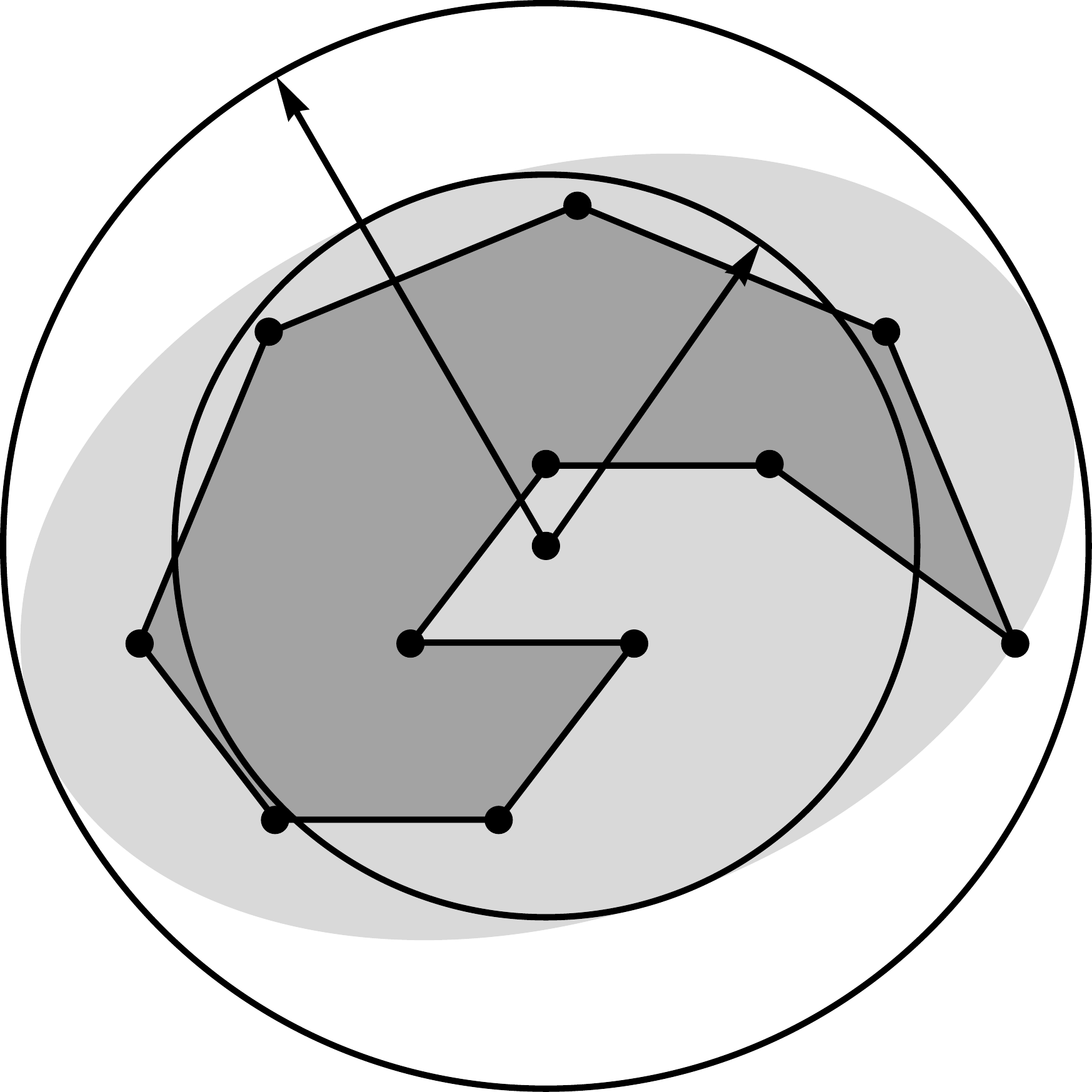}}%
    \put(0.51703827,0.46592345){\color[rgb]{0,0,0}\makebox(0,0)[lt]{\lineheight{1.25}\smash{\begin{tabular}[t]{l}$\bX_c^e$\end{tabular}}}}%
    \put(0.62778703,0.62991683){\color[rgb]{0,0,0}\makebox(0,0)[lt]{\lineheight{1.25}\smash{\begin{tabular}[t]{l}$R_i$\end{tabular}}}}%
    \put(0.33,0.62139767){\color[rgb]{0,0,0}\makebox(0,0)[lt]{\lineheight{1.25}\smash{\begin{tabular}[t]{l}$R_o$\end{tabular}}}}%
  \end{picture}%
\endgroup%

	\caption{\normalsize Sample element enclosed by minimal area ellipse \label{fig:Stab_Geometry}}
\end{figure}
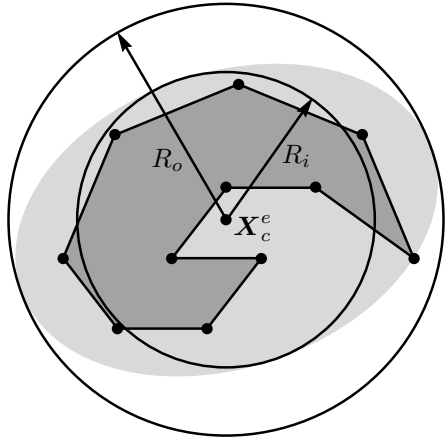

The $\hat{\mu}$ term is scaled by a geometric factor $\beta$, where $\beta=\sqrt{AR}$ with $AR$ the aspect ratio of the element, and an (optional) incompressibility factor $\alpha$. To compute the aspect ratio we make use of a minimal area ellipse enclosing the element. An example of such an ellipse is shown in light grey in Figure \ref{fig:Stab_Geometry}. The centroid of the ellipse is denoted by $\bX_{c}^{e}$ and $R_{i}$ and $R_{o}$ denote the inner and outer radii of the ellipse which are defined by the minor and major radii respectively. The aspect ratio is then computed by $AR=R_{o}/R_{i}$.

The incompressibility factor $\alpha$ is intended to increase the stabilization energy for nearly-incompressible materials and is based on $\TEfive$. However, it is normalised by the Young's modulus such that $\alpha=\TEfive/E_{y}$. The modified $\hat{\mu}$ parameter is the computed from
\begin{equation}
\hat{\mu}=\beta(1+\alpha \beta)\mu \, .
\end{equation}
The merit of the $\alpha$-term is most evident in cases involving complex element geometries and severe deformations and can be set to $0$ for simpler problems. This statement is expanded on in Section \ref{sec:Results}.

We note that it is possible to compute effective values of the Young's modulus and Poisson's ratio from the linearization of any isotropic material model, making the proposed stabilization strain energy widely applicable. 

\section{Numerical Results}
\label{sec:Results}
In this section we present three sample problems, each of which is solved using a different material model, to demonstrate the behaviour of the VEM and the proposed stabilization strain energy function. For each problem we consider the convergence of the displacement at a chosen point, as well as an $\mathcal{H}_{1}$-like error defined by
\begin{align}
||\tilde{\bu}-\bu_{h}||_{1}=&\left[ \int_{\Omega}\left[ |\tilde{\bu}-\bu_{h}|^{2} + |\nabla \tilde{\bu} - \nabla \bu_{h}|^{2} \right] \, d\Omega \right]^{0.5} \nonumber \\
\begin{split}
\approx & \left[ \sum_{i=1}^{n_{el}}\frac{|E_{i}|}{n_{\rm v}^{i}} \sum_{j=1}^{n_{\rm v}} \Bigl[ \left( \tilde{\bu}(\bX_{j}) - \bu_{h}^{i}(\bX_{j}) \right) \cdot \left( \tilde{\bu}(\bX_{j}) -  \bu_{h}^{i}(\bX_{j}) \right) \Bigr. \right. \\
&\Biggl. \Bigl. + \left( \nabla \tilde{\bu}(\bX_{j}) - \Pi \bu_{h}^{i}(\bX_{j}) \right) : \left( \nabla \tilde{\bu}(\bX_{j}) - \Pi \bu_{h}^{i}(\bX_{j}) \right) \Bigr]  \Biggr]^{0.5} \, .
\end{split} 
\label{eqn:H1}
\end{align}
In (\ref{eqn:H1}) $\tilde{\bu}$ denotes a reference solution generated from a mesh of $2^8 \times 2^8$ biquadratic Q2 finite elements, and $\bu_{h}^{i}(\bX_{j})$ denotes the displacement field of element $E_{i}$ evaluated at vertex $\bX_{j}$. Further, $\nabla \tilde{\bu}$ denotes the gradient of the reference solution and  $\Pi \bu_{h}^{i}$ the gradient of $\bu_{h}^{i}$ computed via the projection operator.

\begin{figure}[htp!]
	\begin{subfigure}{.25\textwidth}
		\centering
		\includegraphics[width=0.8\linewidth]{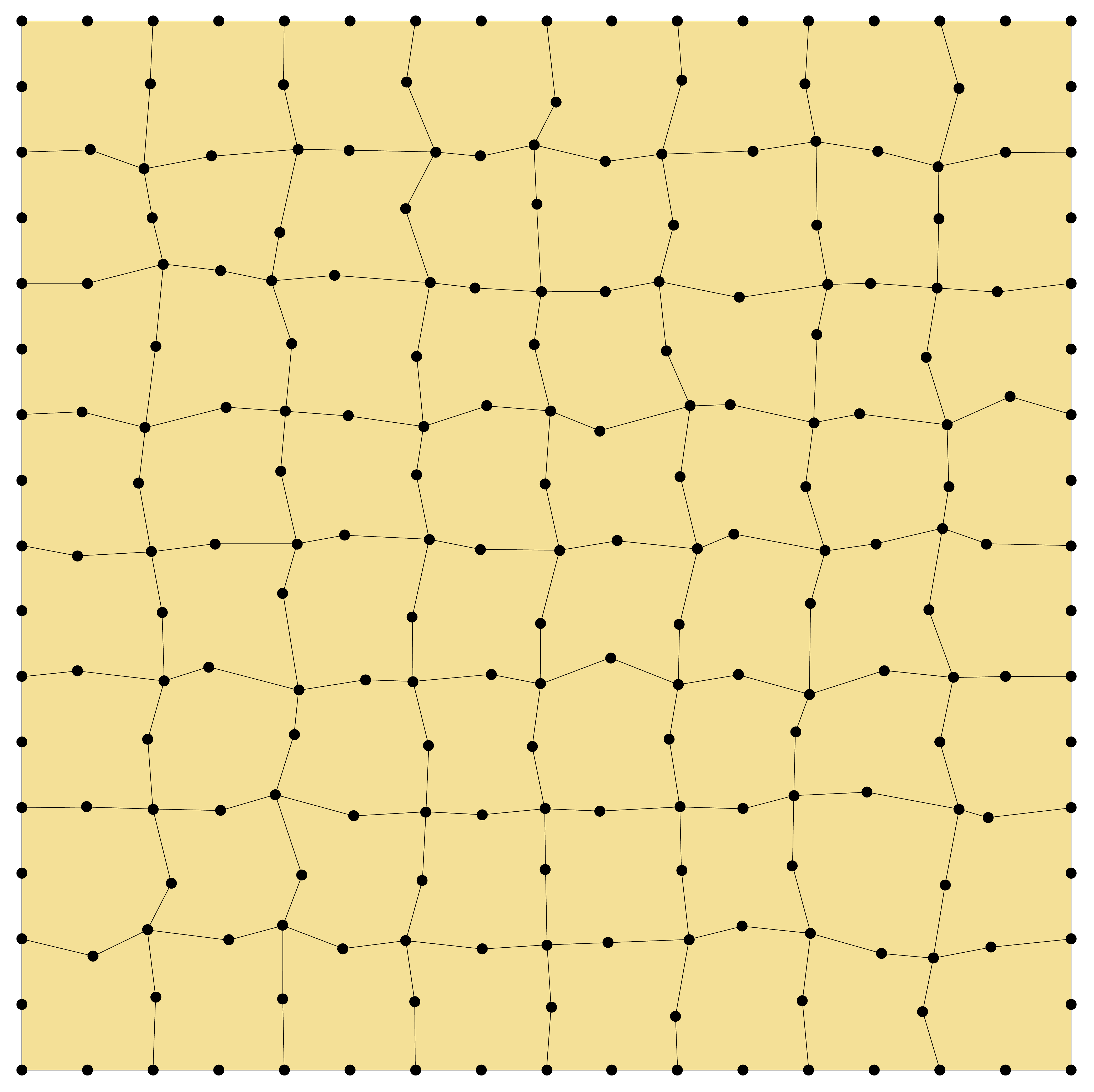}
		\vspace*{-1mm}
		\caption{}
	\end{subfigure}%
	\begin{subfigure}{.25\textwidth}
		\centering
		\includegraphics[width=0.8\linewidth]{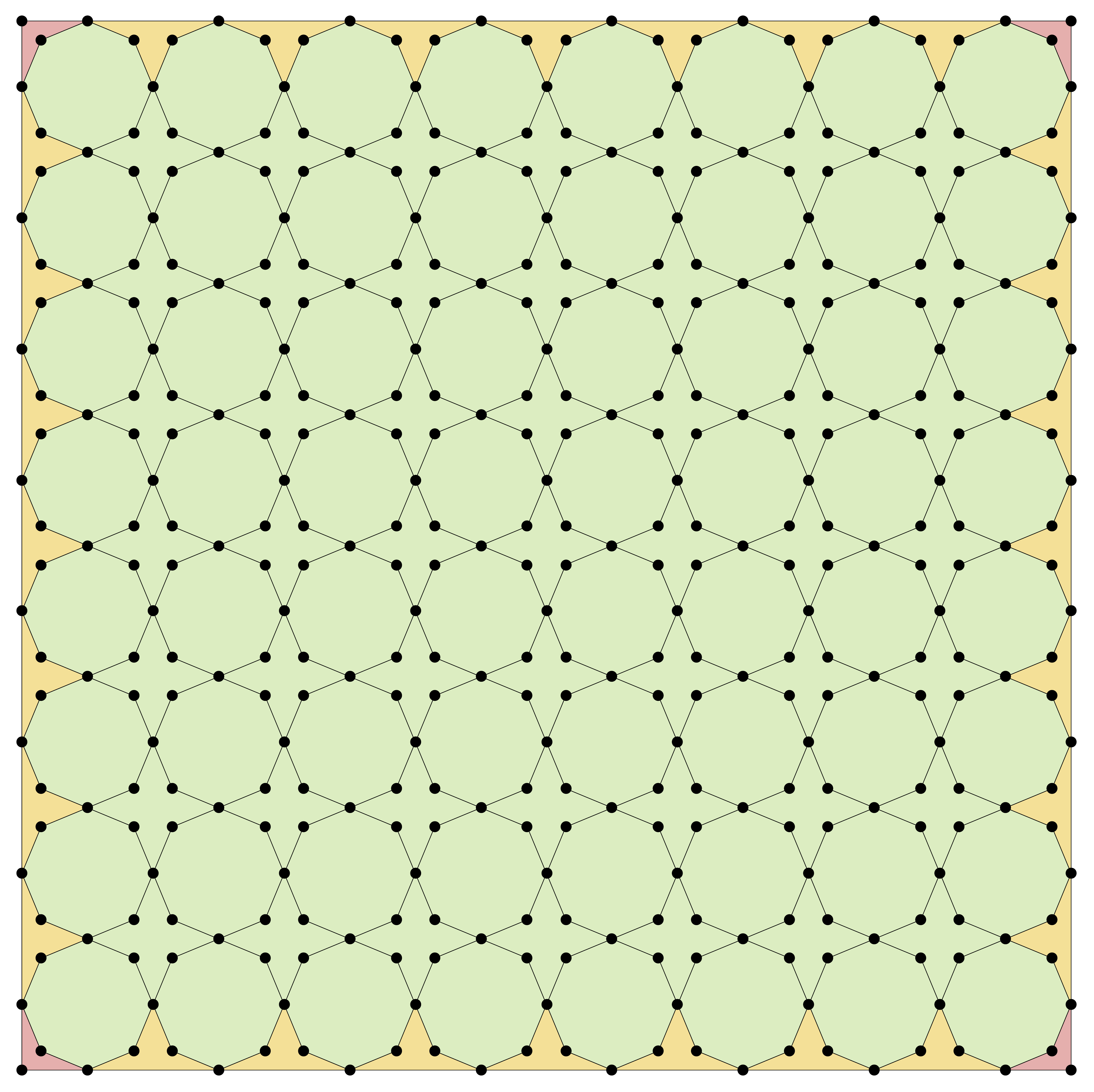}
		\vspace*{-1mm}
		\caption{}
	\end{subfigure}%
	\begin{subfigure}{.25\textwidth}
		\centering
		\includegraphics[width=0.8\linewidth]{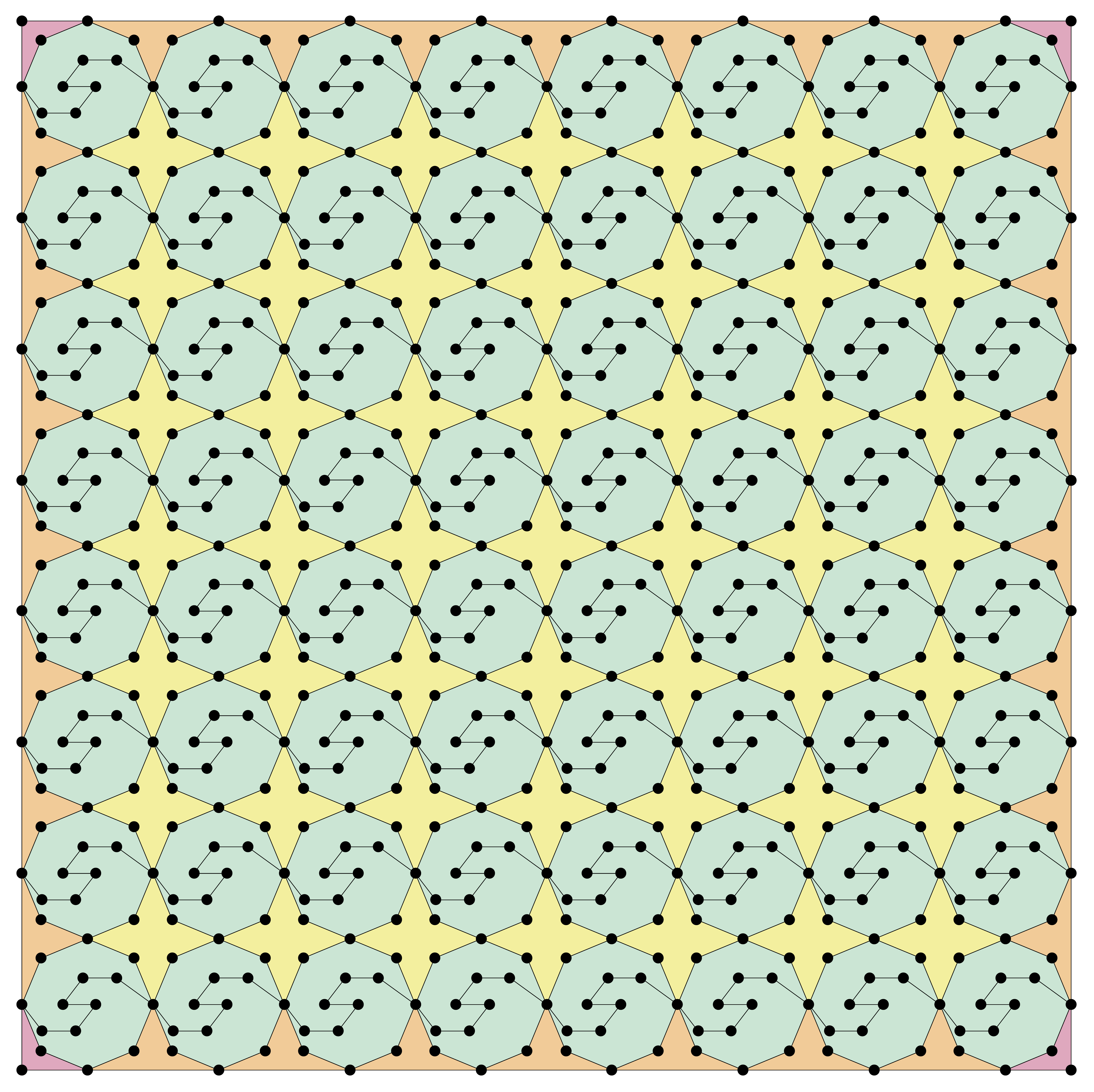}
		\vspace*{-1mm}
		\caption{}
	\end{subfigure}%
	\begin{subfigure}{.25\textwidth}
		\centering
		\includegraphics[width=0.8\linewidth]{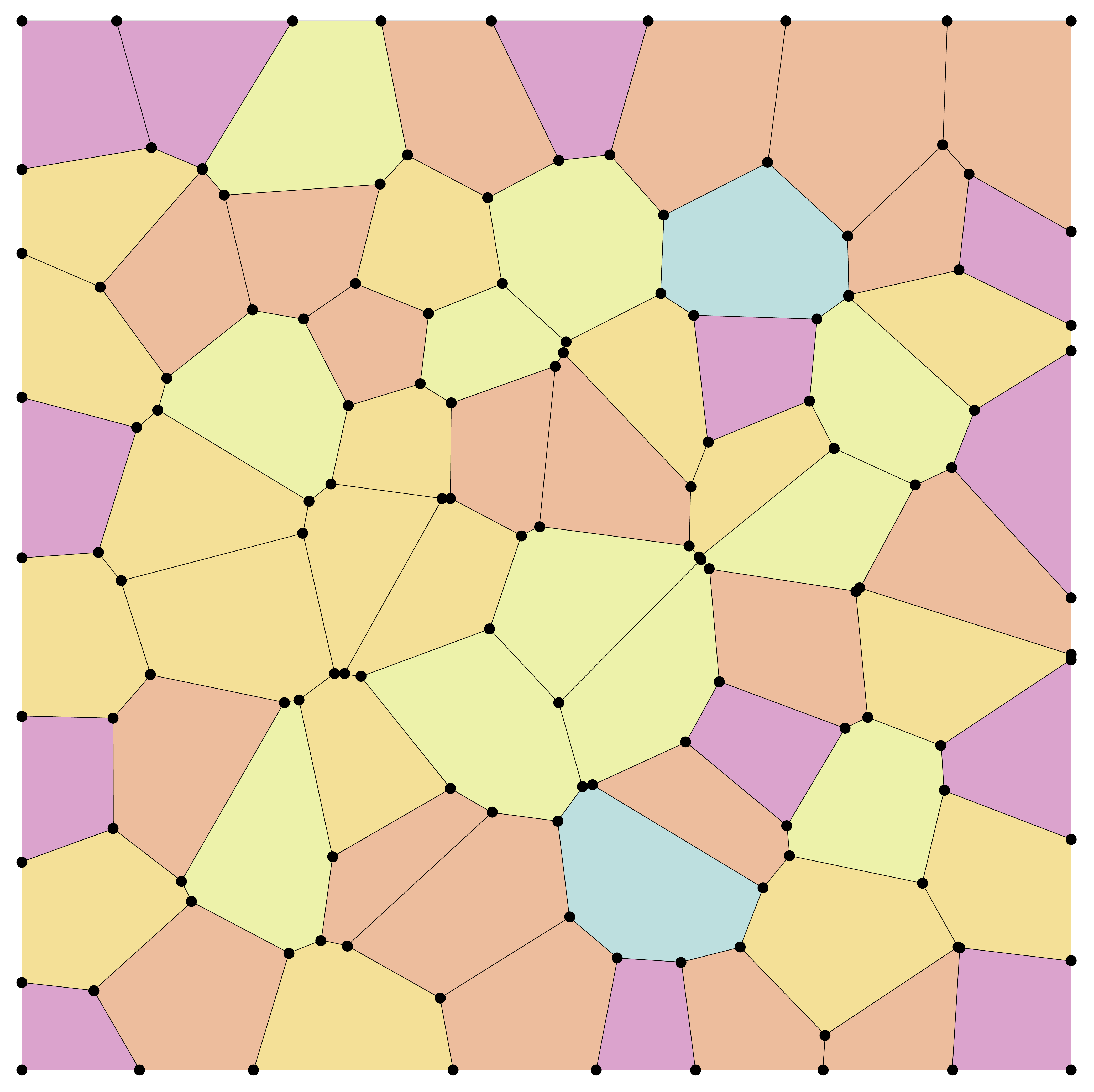}
		\vspace*{-1mm}
		\caption{}
	\end{subfigure}
\vspace*{-3mm}
	\caption{Examples of the various mesh types considered; (a) distorted 8-noded quadrilateral, (b) sun and star, (c) interlocking sun and star, and (d) Voronoi meshes for $N=3$.  \label{fig:MeshTypes}}
\end{figure}

We consider a variety of mesh types with varying degrees of irregularity and concavity; these are depicted in Figure \ref{fig:MeshTypes}. Figure \ref{fig:MeshTypes}(a) shows a distorted 8-noded quadrilateral mesh, denoted by DQ2S. Figure \ref{fig:MeshTypes}(b) shows a mesh comprising convex sun and concave star elements, denoted by S\&S. Figure \ref{fig:MeshTypes}(c) shows an adapted S\&S mesh in which the sun elements are made from two interlocking elements; this is denoted by IS\&S. Finally, Figure \ref{fig:MeshTypes}(d) depicts a classical Voronoi mesh, denoted by VRN. In addition to the meshes shown we consider a structured quadrilateral mesh, denoted by SQ1, and we include the standard Q2 finite element as a reference formulation. For SQ1, DQ2S, VRN and Q2 meshes the number of elements in a mesh is given by $n_{el}=2^{2N}$, where $N$ denotes mesh refinement level, while for S\&S and IS\&S meshes the number of sun elements is given by $n_{el}=2^{2N}$. The meshes depicted in Figure \ref{fig:MeshTypes} correspond to $N=3$. 

In the examples that follow $\alpha=\TEfive/E_{y}$ unless otherwise stated.

\begin{figure}[ht!]
	\centering
	\begin{subfigure}{0.5\textwidth}
		\centering
		\def\svgwidth{0.7\textwidth}
\begingroup%
  \makeatletter%
  \providecommand\color[2][]{%
    \errmessage{(Inkscape) Color is used for the text in Inkscape, but the package 'color.sty' is not loaded}%
    \renewcommand\color[2][]{}%
  }%
  \providecommand\transparent[1]{%
    \errmessage{(Inkscape) Transparency is used (non-zero) for the text in Inkscape, but the package 'transparent.sty' is not loaded}%
    \renewcommand\transparent[1]{}%
  }%
  \providecommand\rotatebox[2]{#2}%
  \newcommand*\fsize{\dimexpr\f@size pt\relax}%
  \newcommand*\lineheight[1]{\fontsize{\fsize}{#1\fsize}\selectfont}%
  \ifx\svgwidth\undefined%
    \setlength{\unitlength}{346.69871612bp}%
    \ifx\svgscale\undefined%
      \relax%
    \else%
      \setlength{\unitlength}{\unitlength * \real{\svgscale}}%
    \fi%
  \else%
    \setlength{\unitlength}{\svgwidth}%
  \fi%
  \global\let\svgwidth\undefined%
  \global\let\svgscale\undefined%
  \makeatother%
  \begin{picture}(1,1.13199224)%
    \lineheight{1}%
    \setlength\tabcolsep{0pt}%
    \put(0,0){\includegraphics[width=\unitlength,page=1]{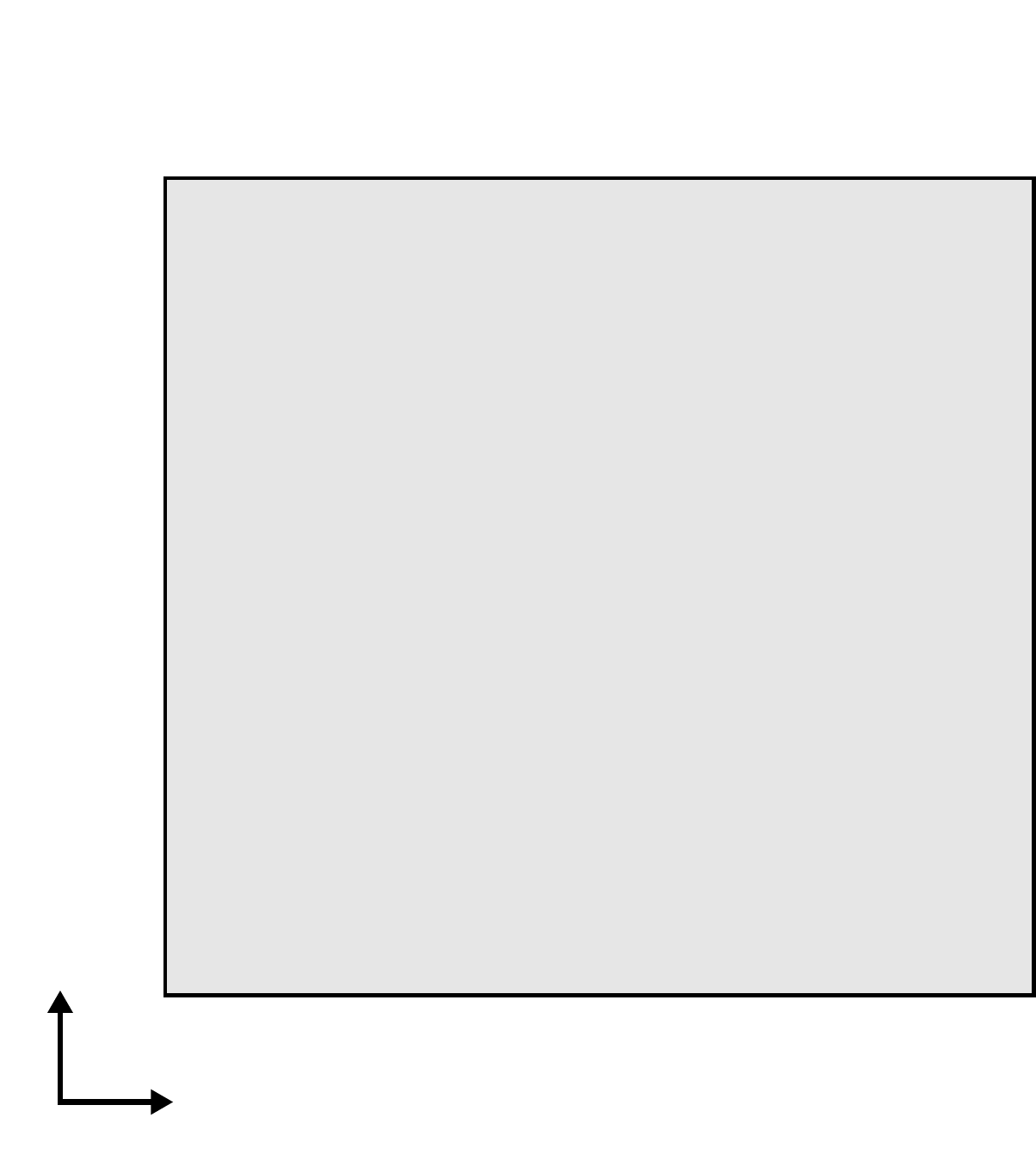}}%
    \put(0.124128,0){\color[rgb]{0,0,0}\makebox(0,0)[lt]{\lineheight{1.25}\smash{\begin{tabular}[t]{l}$x$\end{tabular}}}}%
    \put(-0.00257732,0.14061192){\color[rgb]{0,0,0}\makebox(0,0)[lt]{\lineheight{1.25}\smash{\begin{tabular}[t]{l}$y$\end{tabular}}}}%
    \put(0.54,1.065){\color[rgb]{0,0,0}\makebox(0,0)[lt]{\lineheight{1.25}\smash{\begin{tabular}[t]{l}$q_{\rm s}$\end{tabular}}}}%
    \put(0,0){\includegraphics[width=\unitlength,page=2]{isotropic_simple_shear.pdf}}%
  \end{picture}%
\endgroup%

		\vspace*{-6mm}
		\caption{}
	\end{subfigure}%
	\begin{subfigure}{0.5\textwidth}
		\centering
		\includegraphics[width=0.9\textwidth]{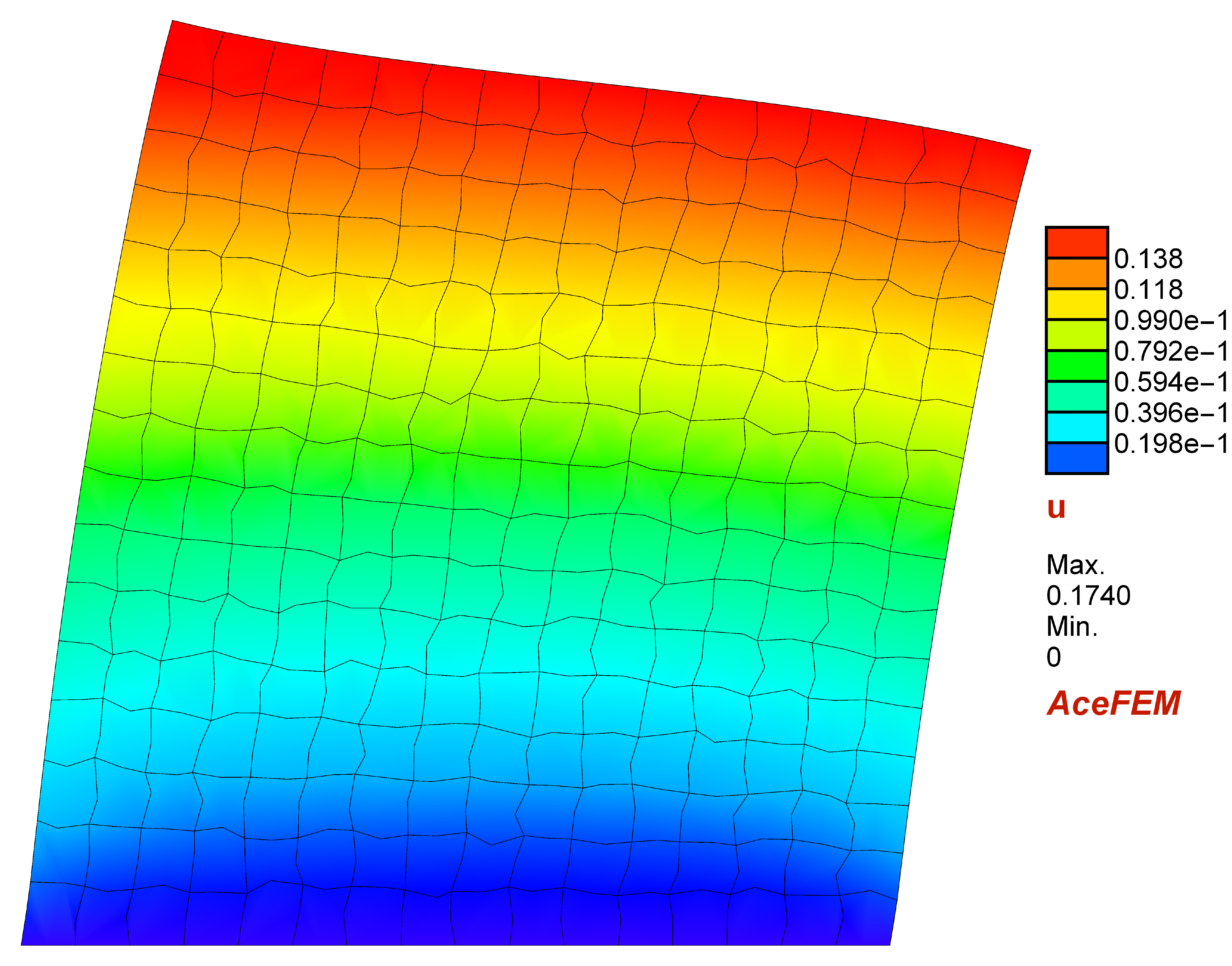}
		\caption{}
	\end{subfigure}
	\caption{\normalsize Simple shear test (a) problem geometry and (b) sample deformed configuration for a neo-Hookean material model with $N=4$ and $\nu =0.3$ for a DQ2S mesh. \label{fig:SimpleShear}} \vspace*{-2mm}
\end{figure}

\begin{figure}[ht!]
	\centering
	\begin{subfigure}{0.5\textwidth}
		\centering
		\def\svgwidth{\textwidth}
		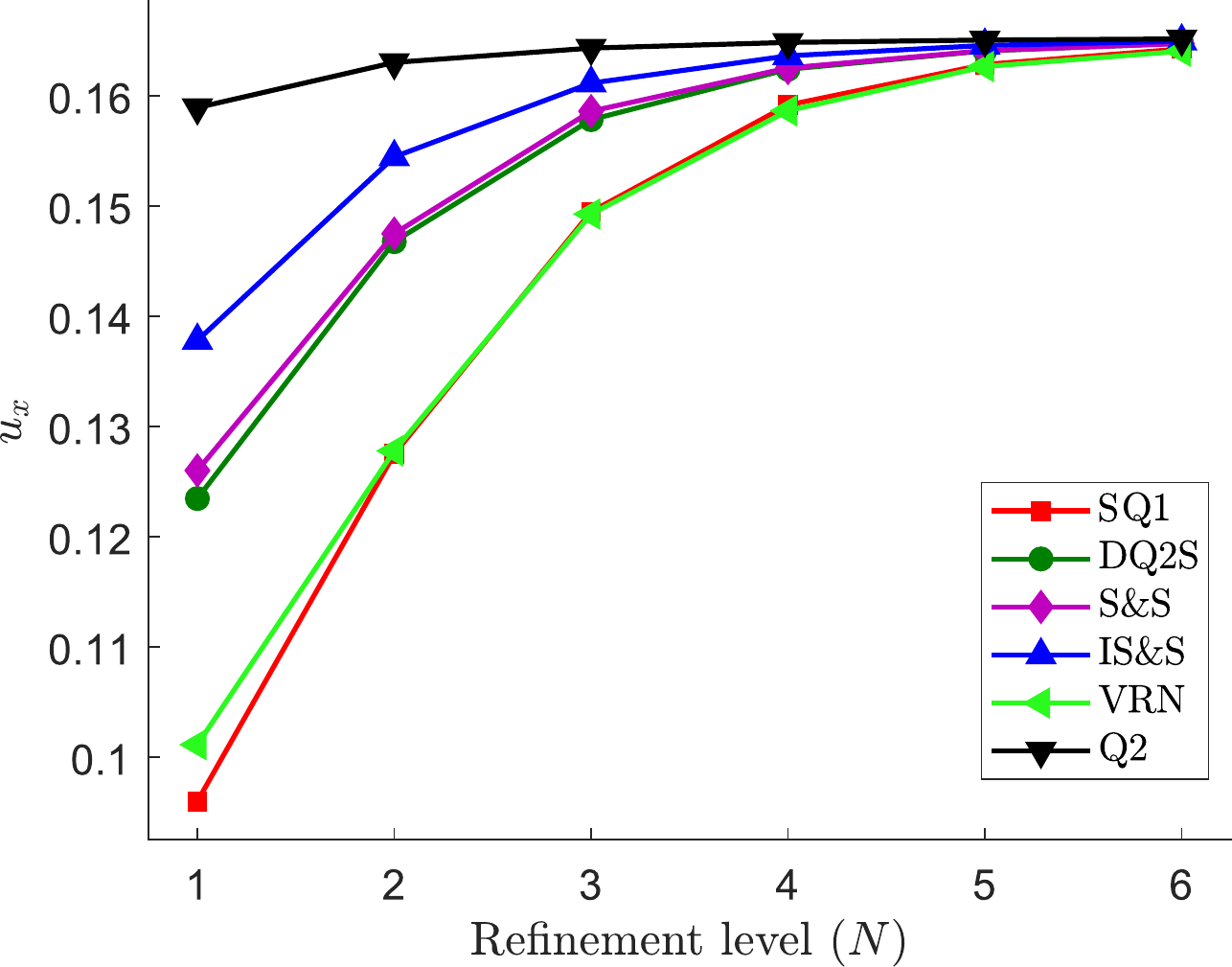
		\vspace*{-6mm}
		\caption{}
	\end{subfigure}%
	\begin{subfigure}{0.5\textwidth}
		\centering
		\includegraphics[width=\textwidth]{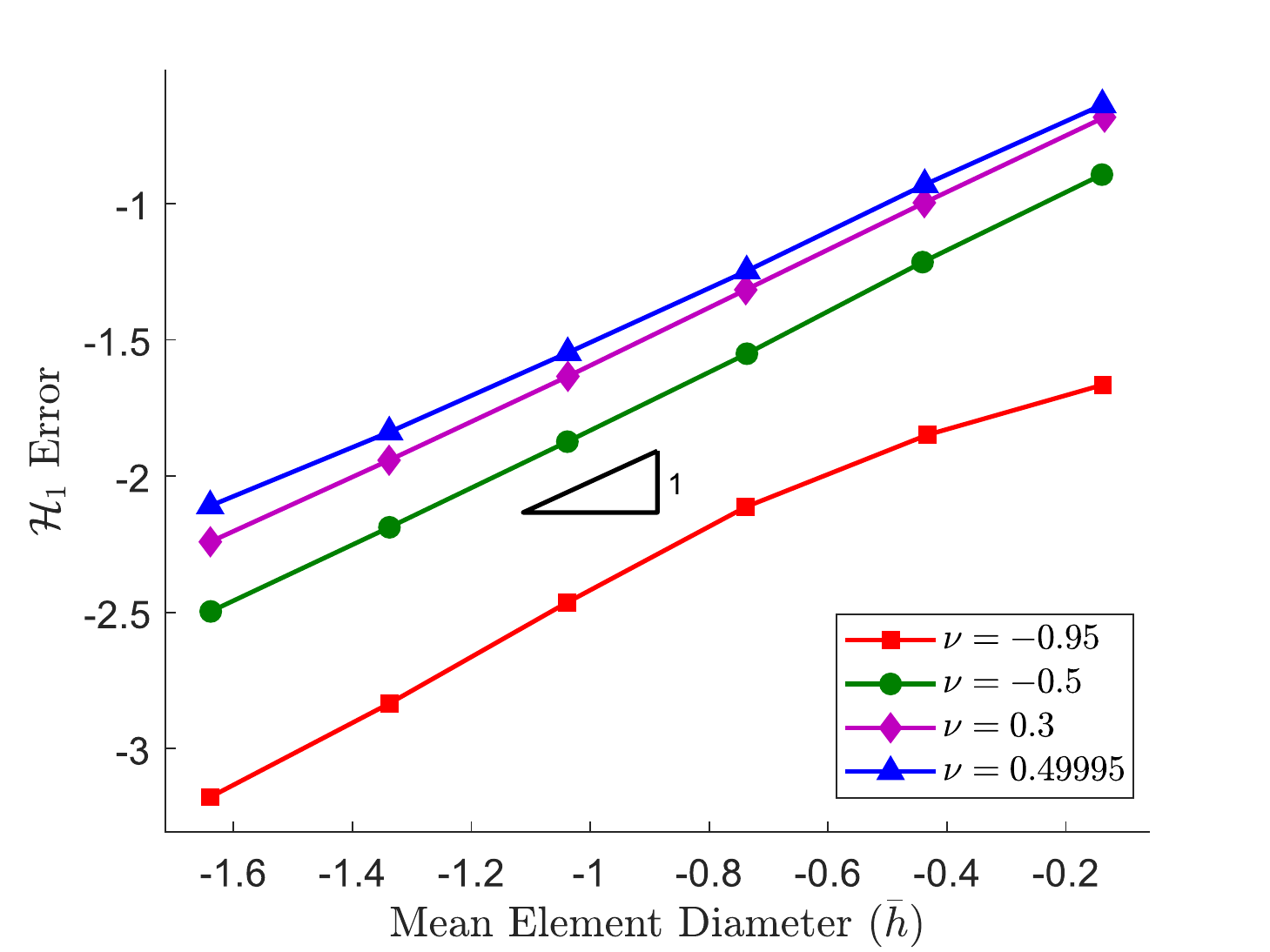}
		\vspace*{-6mm}
		\caption{}
	\end{subfigure}
	\caption{\normalsize Simple shear test (a) tip deflection vs mesh refinement for $\nu=0.3$ and (b) loglog $\mathcal{H}_{1}$ error vs mean element diameter $\bar{h}$ with DQ2S meshes, for a neo-Hookean material model. \label{fig:SimpleShearConvergence}} \vspace*{-2mm} 
\end{figure}

\paragraph*{Simple shear} This problem consists of a unit square domain, fully constrained along the bottom, subjected to a uniformly distributed horizontal load $q_{\rm s}=5~\frac{\rm N}{\rm m}$ along the top, as shown in Figure \ref{fig:SimpleShear}(a). For this problem we consider a neo-Hookean material model with the strain energy function 
\begin{align}
\Psi_{\rm NH}&=\frac{\mu}{2}(I_{C} - 3 - 2 \ln (J)) + \frac{\lambda}{2}\ln (J)^{2} \, . \label{eqn:NeoHookComp}
\end{align} 
We choose a value of Young's modulus of $E_{y}=200~\rm Pa$. Figure \ref{fig:SimpleShear}(b) shows the deformed configuration of the simple shear problem with a distorted quadrilateral DQ2S mesh for $N=4$ and a value of Poisson's ratio of $\nu=0.3$.
\\
 \\
Figure \ref{fig:SimpleShearConvergence}(a) shows a graph of the horizontal displacement of the upper right hand corner of the body vs mesh refinement level with a Poisson's ratio of $\nu =0.3$. We note smooth monotonic convergence of all formulations and mesh types. Figure \ref{fig:SimpleShearConvergence}(b) shows a plot of $\mathcal{H}_{1}$ error vs mean element diameter $\bar{h}$ on a loglog scale for the simple shear problem, with DQ2S meshes, for a variety of choices of Poisson's ratio. We note, for most choices of Poisson's ratio, a near-optimal convergence rate with a gradient of approximately 1. However, for $\nu=-0.95$ the coarse mesh convergence behaviour is non-optimal; this is likely due to the `over-stabilization' discussed in Section \ref{subsec:Stab}.
\newpage

\begin{figure}[htp!]
	\centering
	\begin{subfigure}{.5\textwidth}
		\centering
		\def\svgwidth{0.8\textwidth}
\begingroup%
  \makeatletter%
  \providecommand\color[2][]{%
    \errmessage{(Inkscape) Color is used for the text in Inkscape, but the package 'color.sty' is not loaded}%
    \renewcommand\color[2][]{}%
  }%
  \providecommand\transparent[1]{%
    \errmessage{(Inkscape) Transparency is used (non-zero) for the text in Inkscape, but the package 'transparent.sty' is not loaded}%
    \renewcommand\transparent[1]{}%
  }%
  \providecommand\rotatebox[2]{#2}%
  \newcommand*\fsize{\dimexpr\f@size pt\relax}%
  \newcommand*\lineheight[1]{\fontsize{\fsize}{#1\fsize}\selectfont}%
  \ifx\svgwidth\undefined%
    \setlength{\unitlength}{560.99551413bp}%
    \ifx\svgscale\undefined%
      \relax%
    \else%
      \setlength{\unitlength}{\unitlength * \real{\svgscale}}%
    \fi%
  \else%
    \setlength{\unitlength}{\svgwidth}%
  \fi%
  \global\let\svgwidth\undefined%
  \global\let\svgscale\undefined%
  \makeatother%
  \begin{picture}(1,1.02857725)%
    \lineheight{1}%
    \setlength\tabcolsep{0pt}%
    \put(0,0){\includegraphics[width=\unitlength,page=1]{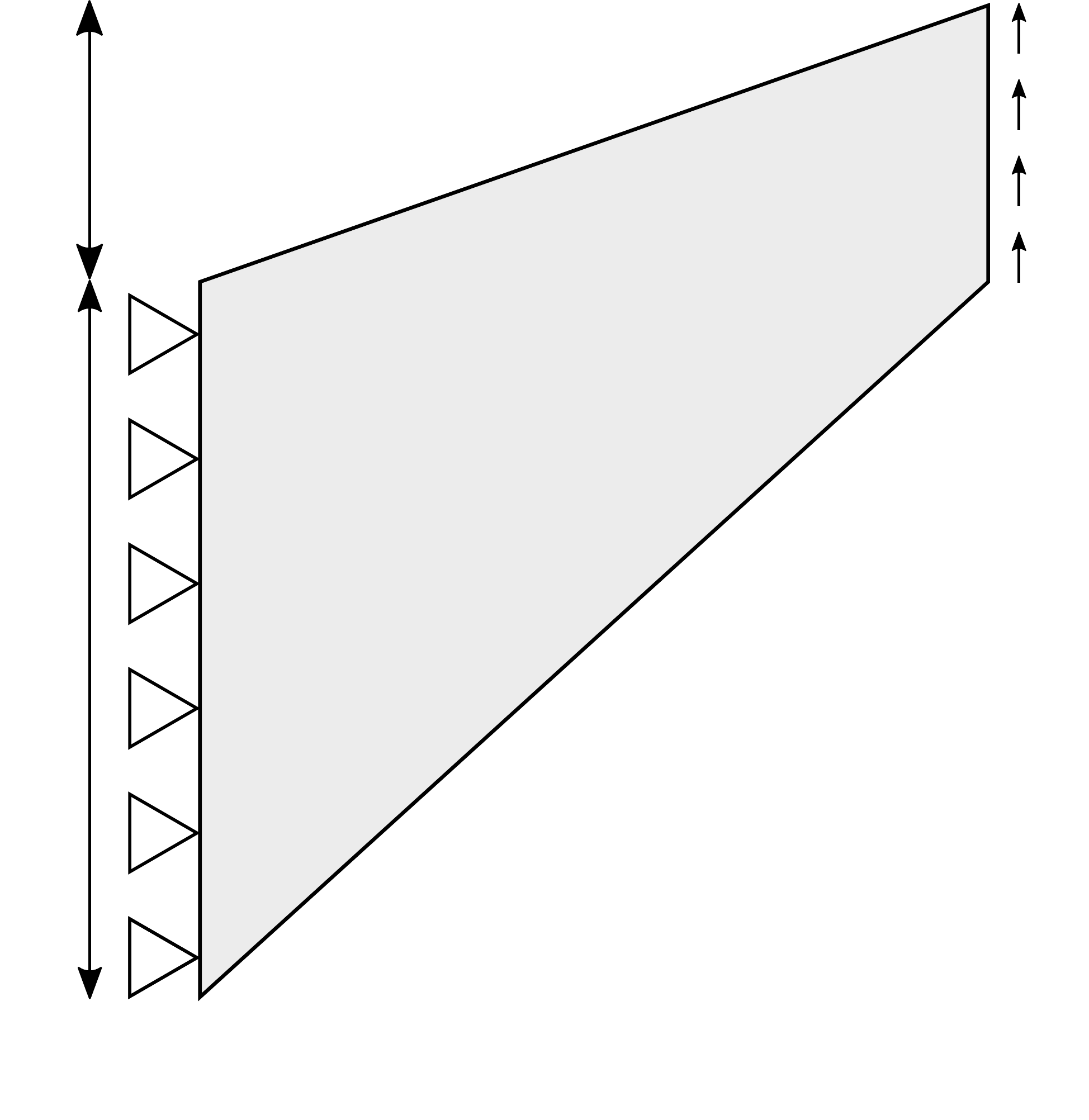}}%
    \put(0.94952255,0.88426986){\color[rgb]{0,0,0}\makebox(0,0)[lt]{\lineheight{1.25}\smash{\begin{tabular}[t]{l}$q_{\rm c}$\end{tabular}}}}%
    \put(0,0){\includegraphics[width=\unitlength,page=2]{cook.pdf}}%
    \put(0.51789184,0.11115442){\color[rgb]{0,0,0}\makebox(0,0)[lt]{\lineheight{1.25}\smash{\begin{tabular}[t]{l}$60$\end{tabular}}}}%
    \put(0.00604667,0.42208133){\color[rgb]{0,0,0}\makebox(0,0)[lt]{\lineheight{1.25}\smash{\begin{tabular}[t]{l}$44$\end{tabular}}}}%
    \put(0.00795657,0.88426986){\color[rgb]{0,0,0}\makebox(0,0)[lt]{\lineheight{1.25}\smash{\begin{tabular}[t]{l}$16$\end{tabular}}}}%
    \put(0,0){\includegraphics[width=\unitlength,page=3]{cook.pdf}}%
    \put(0.13400799,0){\color[rgb]{0,0,0}\makebox(0,0)[lt]{\lineheight{1.25}\smash{\begin{tabular}[t]{l}$x$\end{tabular}}}}%
    \put(-0.0015928,0.13942041){\color[rgb]{0,0,0}\makebox(0,0)[lt]{\lineheight{1.25}\smash{\begin{tabular}[t]{l}$y$\end{tabular}}}}%
  \end{picture}%
\endgroup%

		\vspace*{-6mm}
		\caption{}
	\end{subfigure}%
	\begin{subfigure}{.5\textwidth}
		\centering
		\includegraphics[width=0.7\linewidth]{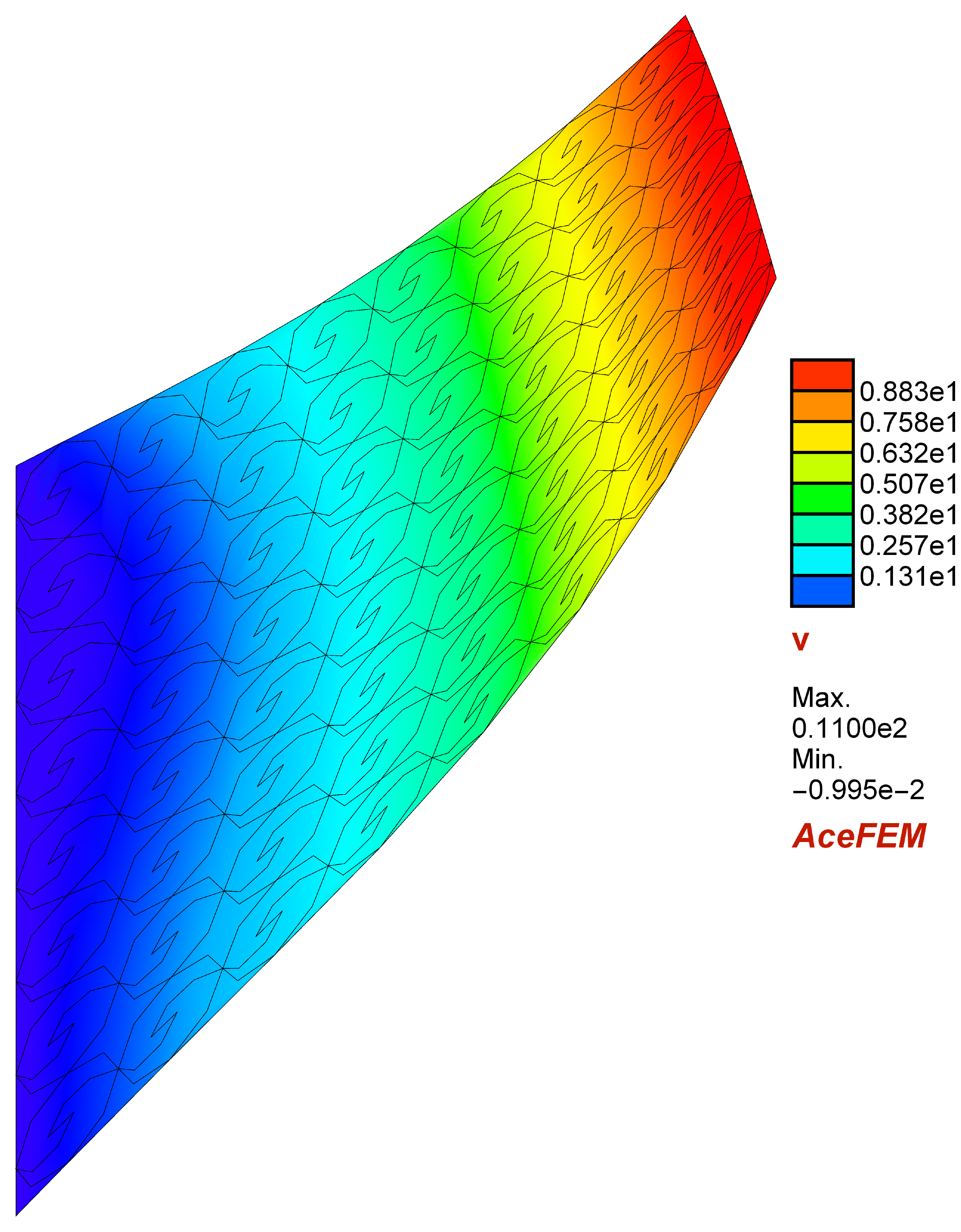}
		\vspace*{-6mm}
		\caption{}
	\end{subfigure}
	\caption{Cook's membrane problem (a) geometry and (b) sample deformed configuration for a Mooney-Rivlin material model with $N=3$ and $\nu=0.49995$ for an IS\&S mesh. \label{fig:Cook}}
\end{figure}

\paragraph*{Cook's Membrane} This problem consists of a uniformly tapered panel, fully constrained on its left edge, subject to a uniformly distributed load $q_{\rm c}$ with a total magnitude $Q_{\rm c}=10~N$ on its right edge, as depicted in Figure \ref{fig:Cook}(a). For this problem we consider a Mooney-Rivlin material model with the strain energy function
\begin{equation}
\Psi_{\rm MR} =  C_{10}(\overline{I}_{C} - 3) + C_{01}(\overline{II}_{C} - 3) + \frac{\kappa}{2}\ln (J)^{2} \, . \label{eqn:MooneyRivlinComp}
\end{equation}
Comparison with Hooke's law yields the relations 
\begin{equation}
\kappa = \frac{E_{y}}{3(1-2\nu)} \quad \text{and} \quad \mu = 2(C_{10}+C_{01}) \, ,
\end{equation}
where, for simplicity, we set $C_{10}=rC_{01}$ and choose $r=4$. We choose a value of Young's modulus of $E_{y}=200~\rm Pa$. Figure \ref{fig:Cook}(b) shows the deformed configuration of Cook's membrane problem with an interlocking sun and star IS\&S mesh for $N=3$ with a Poisson's ratio of $\nu =0.49995$.
\\

\begin{figure}[ht!]
	\centering
	\begin{subfigure}{0.5\textwidth}
		\centering
		\def\svgwidth{\textwidth}
		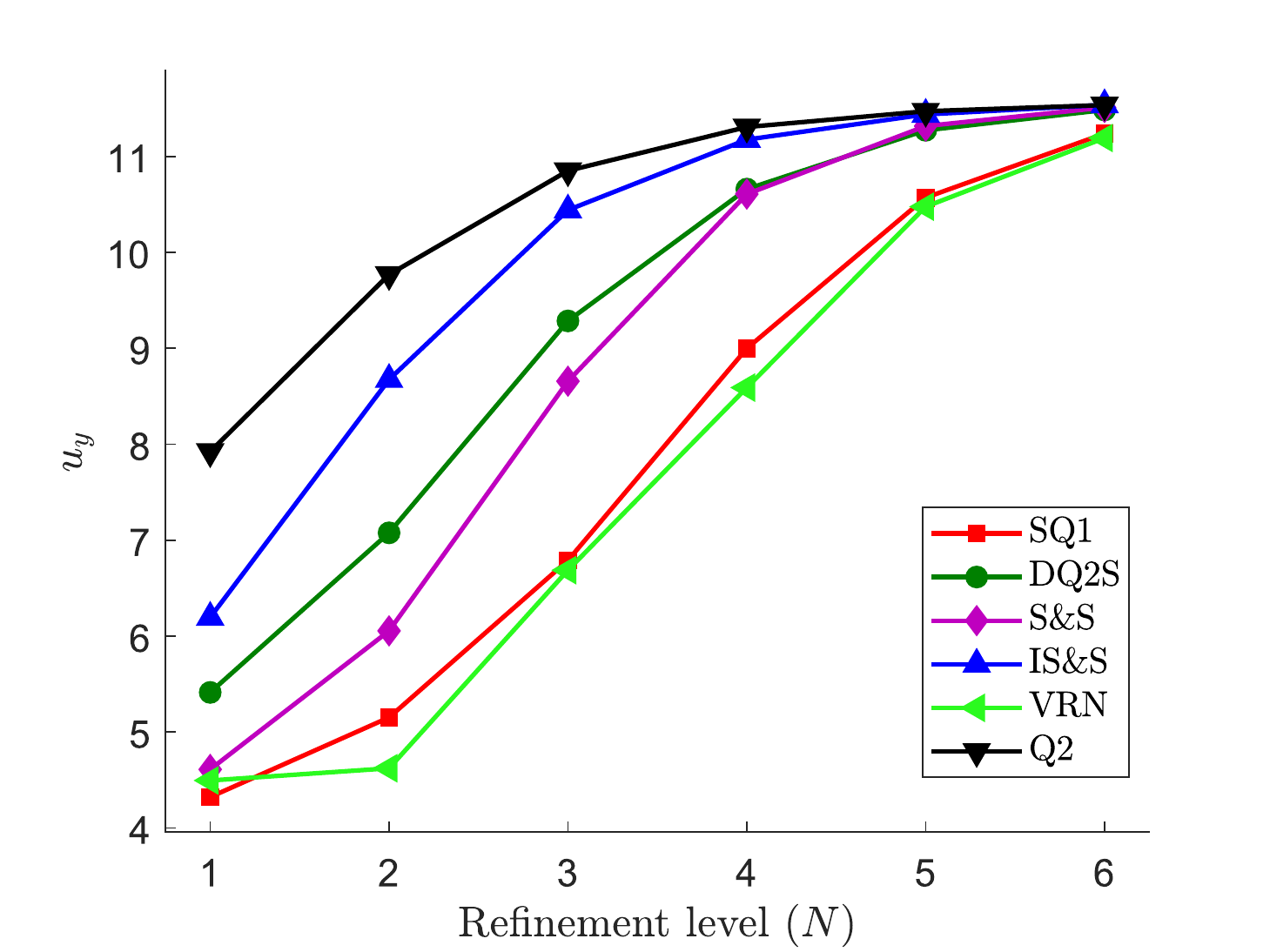
		\vspace*{-5mm}
		\caption{}
	\end{subfigure}%
	\begin{subfigure}{0.5\textwidth}
		\centering
		\def\svgwidth{\textwidth}
		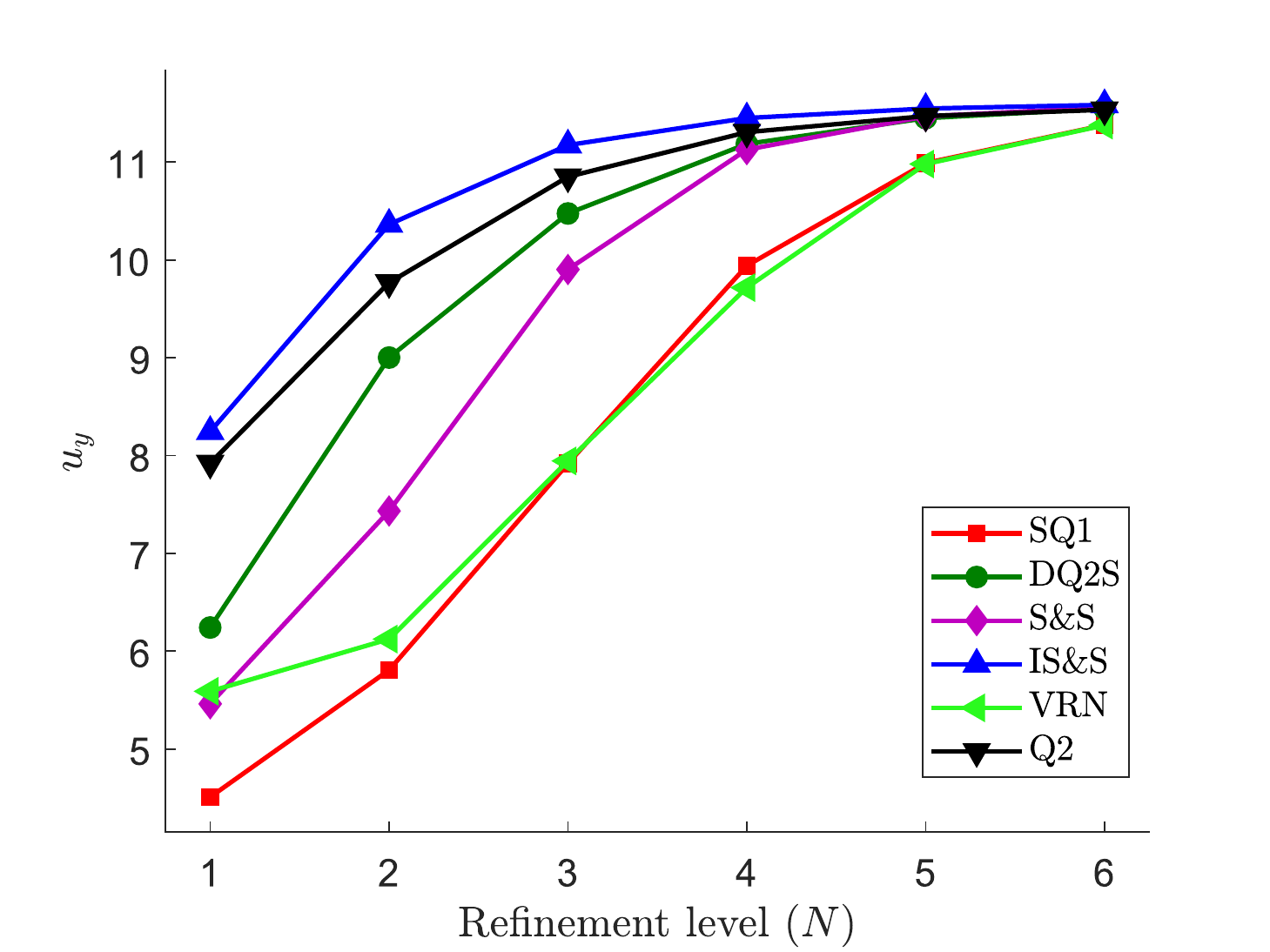
		\vspace*{-5mm}
		\caption{}
	\end{subfigure}
	\caption{\normalsize Cook's membrane problem tip deflection vs mesh refinement for $\nu=0.49995$ (a) with, and (b) without the incompressibility factor for a Mooney-Rivlin material model. \label{fig:CookPointConvergence}} \vspace*{-2mm}
\end{figure}

\begin{figure}[ht!]
	\centering
	\begin{subfigure}{0.5\textwidth}
		\centering
		\includegraphics[width=0.9\textwidth]{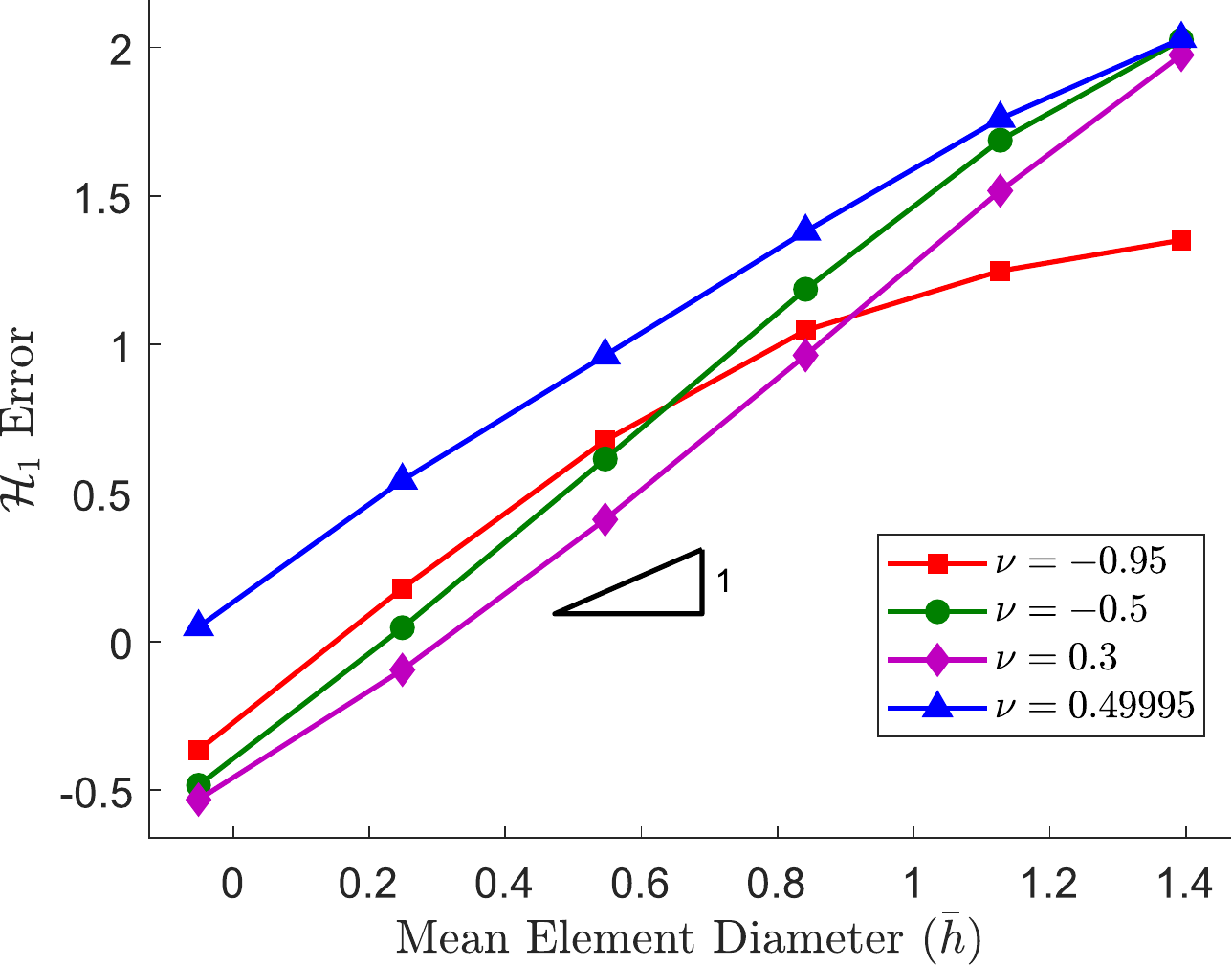}
		\caption{}
	\end{subfigure}%
	\begin{subfigure}{0.5\textwidth}
		\centering
		\includegraphics[width=\textwidth]{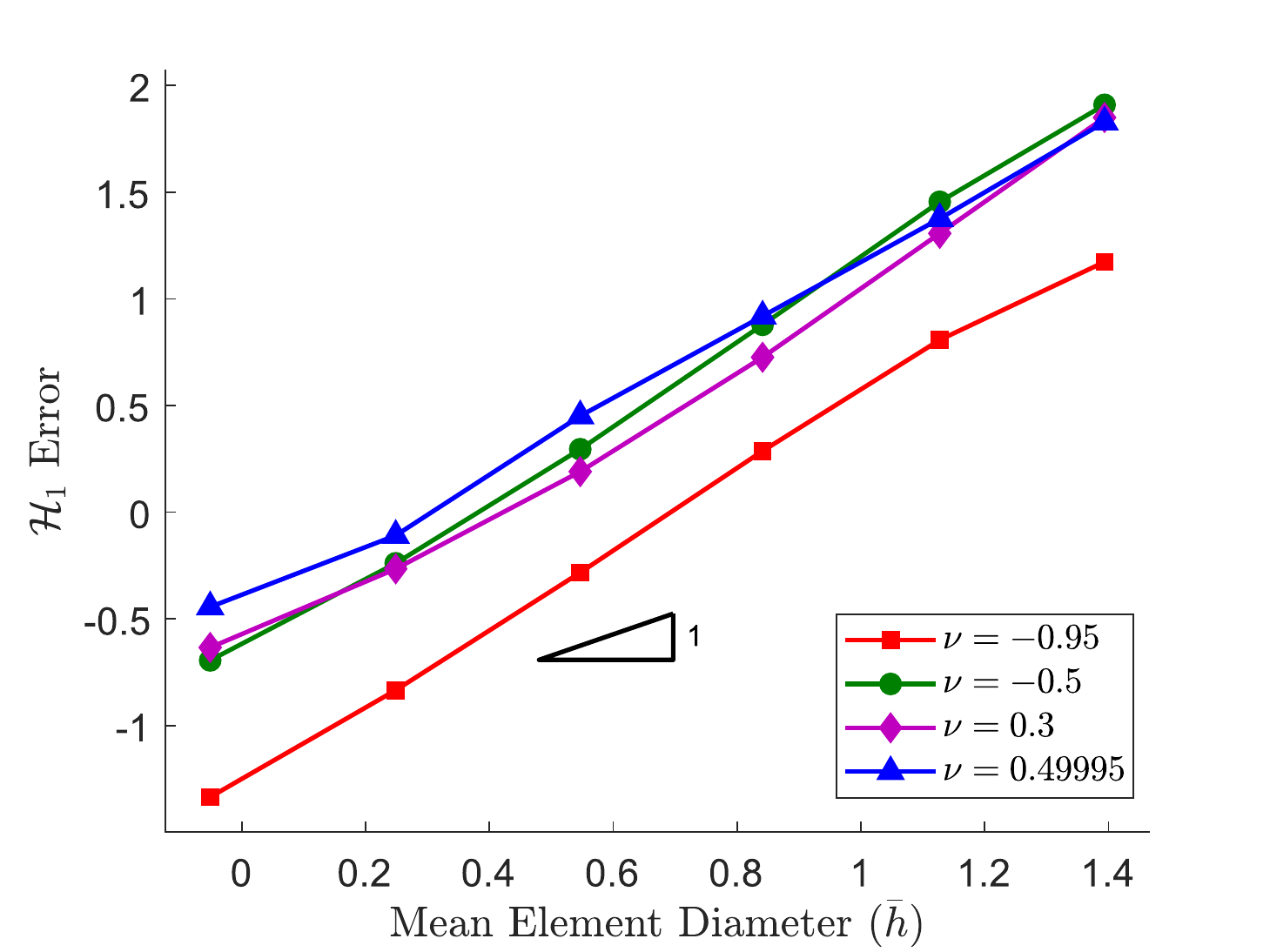}
		\vspace*{-5mm}
		\caption{}
	\end{subfigure}
	\caption{\normalsize Cook's membrane problem loglog $\mathcal{H}_{1}$ error vs mean element diameter with IS\&S meshes, for a Mooney-Rivlin material model. \label{fig:CookH1Convergence}} \vspace*{-2mm}
\end{figure}

\newpage
Figure \ref{fig:CookPointConvergence}(a) shows the vertical displacement of the upper right hand corner of the body vs mesh refinement level with a Poisson's ratio of $\nu =0.49995$. We again note smooth monotonic convergence of all formulations and mesh types. Furthermore, we note that the VEM formulation is locking free in the near-incompressible limit. The SQ1 and VRN meshes show poorer accuracy than the other mesh types as they have significantly fewer degrees of freedom per element. Figure \ref{fig:CookPointConvergence}(b) is similar to Figure \ref{fig:CookPointConvergence}(a), however, we set $\alpha=0$ which `turns off' the incompressibility factor. This softens the stabilization parameters and we note improved accuracy from all VEM formulations. The SQ1 and VRN meshes, however, still exhibit lower accuracy than the other mesh types. 
\\
\newpage
Figure \ref{fig:CookH1Convergence}(a) shows a plot of $\mathcal{H}_{1}$ error vs mean element diameter $\bar{h}$ on a loglog scale for Cook's membrane problem with IS\&S meshes for a variety of choices of Poisson's ratio. Again, Figure \ref{fig:CookH1Convergence}(b) is similar to Figure \ref{fig:CookH1Convergence}(a), however, we set $\alpha=0$. For the case $\alpha=\TEfive/E_{y}$, Figure \ref{fig:CookH1Convergence}(a), we again note suboptimal coarse mesh convergence behaviour for $\nu=-0.95$. For the case $\alpha=0$ this behaviour is significantly improved. However, the fine mesh convergence behaviour for $\nu=0.3$ and $\nu=0.49995$ is less uniform than that in Figure \ref{fig:CookH1Convergence}(a). Interestingly, we observe a greater than expected order of accuracy with a convergence rate of approximately 1.5. 

\begin{figure}[htp!]
	\centering
	\begin{subfigure}{.5\textwidth}
		\centering
		\def\svgwidth{0.8\textwidth}
\begingroup%
  \makeatletter%
  \providecommand\color[2][]{%
    \errmessage{(Inkscape) Color is used for the text in Inkscape, but the package 'color.sty' is not loaded}%
    \renewcommand\color[2][]{}%
  }%
  \providecommand\transparent[1]{%
    \errmessage{(Inkscape) Transparency is used (non-zero) for the text in Inkscape, but the package 'transparent.sty' is not loaded}%
    \renewcommand\transparent[1]{}%
  }%
  \providecommand\rotatebox[2]{#2}%
  \newcommand*\fsize{\dimexpr\f@size pt\relax}%
  \newcommand*\lineheight[1]{\fontsize{\fsize}{#1\fsize}\selectfont}%
  \ifx\svgwidth\undefined%
    \setlength{\unitlength}{487.13162203bp}%
    \ifx\svgscale\undefined%
      \relax%
    \else%
      \setlength{\unitlength}{\unitlength * \real{\svgscale}}%
    \fi%
  \else%
    \setlength{\unitlength}{\svgwidth}%
  \fi%
  \global\let\svgwidth\undefined%
  \global\let\svgscale\undefined%
  \makeatother%
  \begin{picture}(1,0.68928406)%
    \lineheight{1}%
    \setlength\tabcolsep{0pt}%
    \put(0,0){\includegraphics[width=\unitlength,page=1]{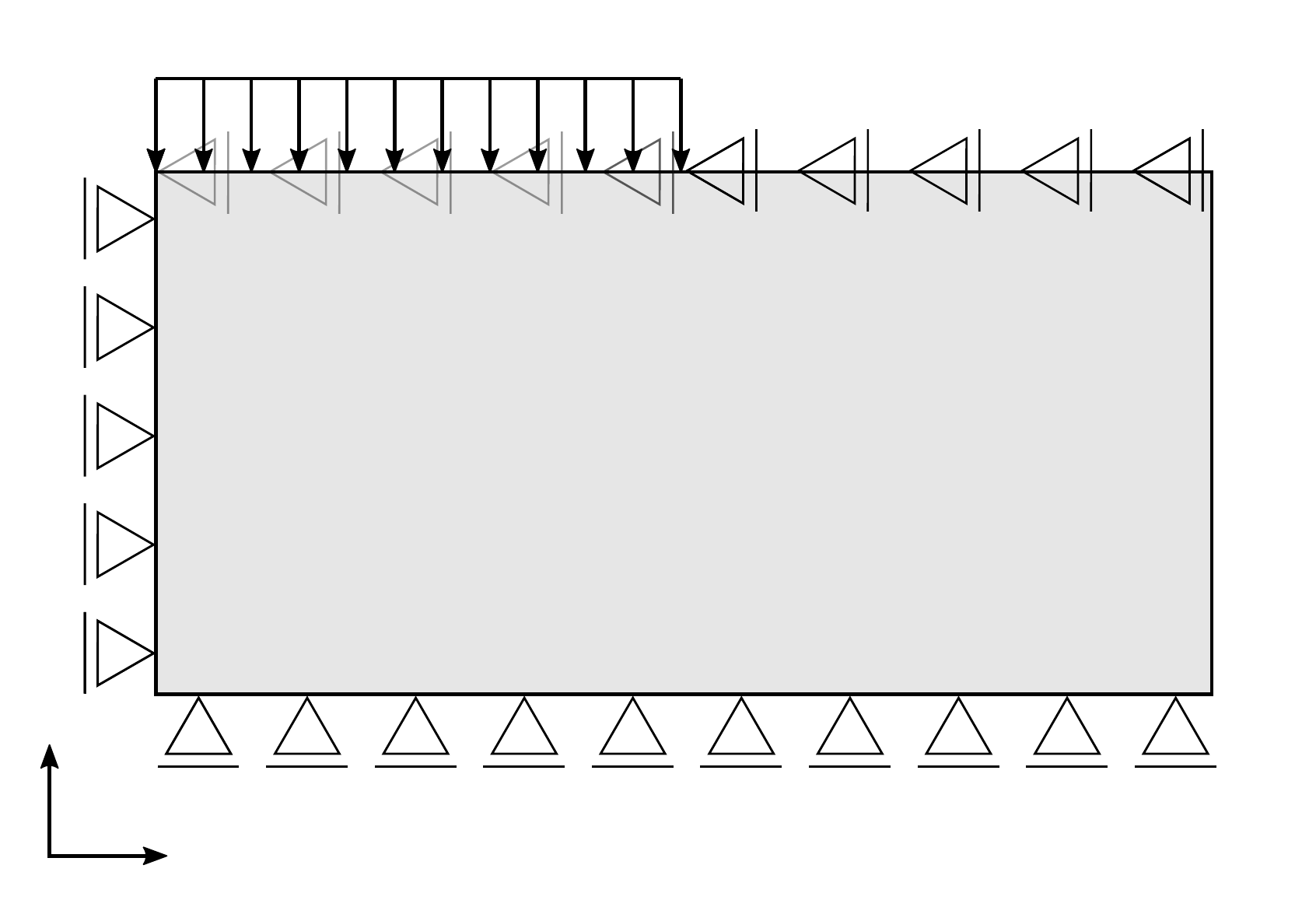}}%
    \put(0.11386607,0){\color[rgb]{0,0,0}\makebox(0,0)[lt]{\lineheight{1.25}\smash{\begin{tabular}[t]{l}$x$\end{tabular}}}}%
    \put(-0.00122288,0.1088679){\color[rgb]{0,0,0}\makebox(0,0)[lt]{\lineheight{1.25}\smash{\begin{tabular}[t]{l}$y$\end{tabular}}}}%
    \put(0.27405743,0.65476278){\color[rgb]{0,0,0}\makebox(0,0)[lt]{\lineheight{1.25}\smash{\begin{tabular}[t]{l}$q_{\rm p}$\end{tabular}}}}%
    \put(0,0){\includegraphics[width=\unitlength,page=2]{punch.pdf}}%
    \put(0.49956956,0.02488401){\color[rgb]{0,0,0}\makebox(0,0)[lt]{\lineheight{1.25}\smash{\begin{tabular}[t]{l}$2$\end{tabular}}}}%
    \put(0,0){\includegraphics[width=\unitlength,page=3]{punch.pdf}}%
    \put(0.9739226,0.34371159){\color[rgb]{0,0,0}\makebox(0,0)[lt]{\lineheight{1.25}\smash{\begin{tabular}[t]{l}$1$\end{tabular}}}}%
  \end{picture}%
\endgroup%

		\vspace*{-2mm}
		\caption{}
	\end{subfigure}%
	\begin{subfigure}{.5\textwidth}
		\centering
		\includegraphics[width=\linewidth]{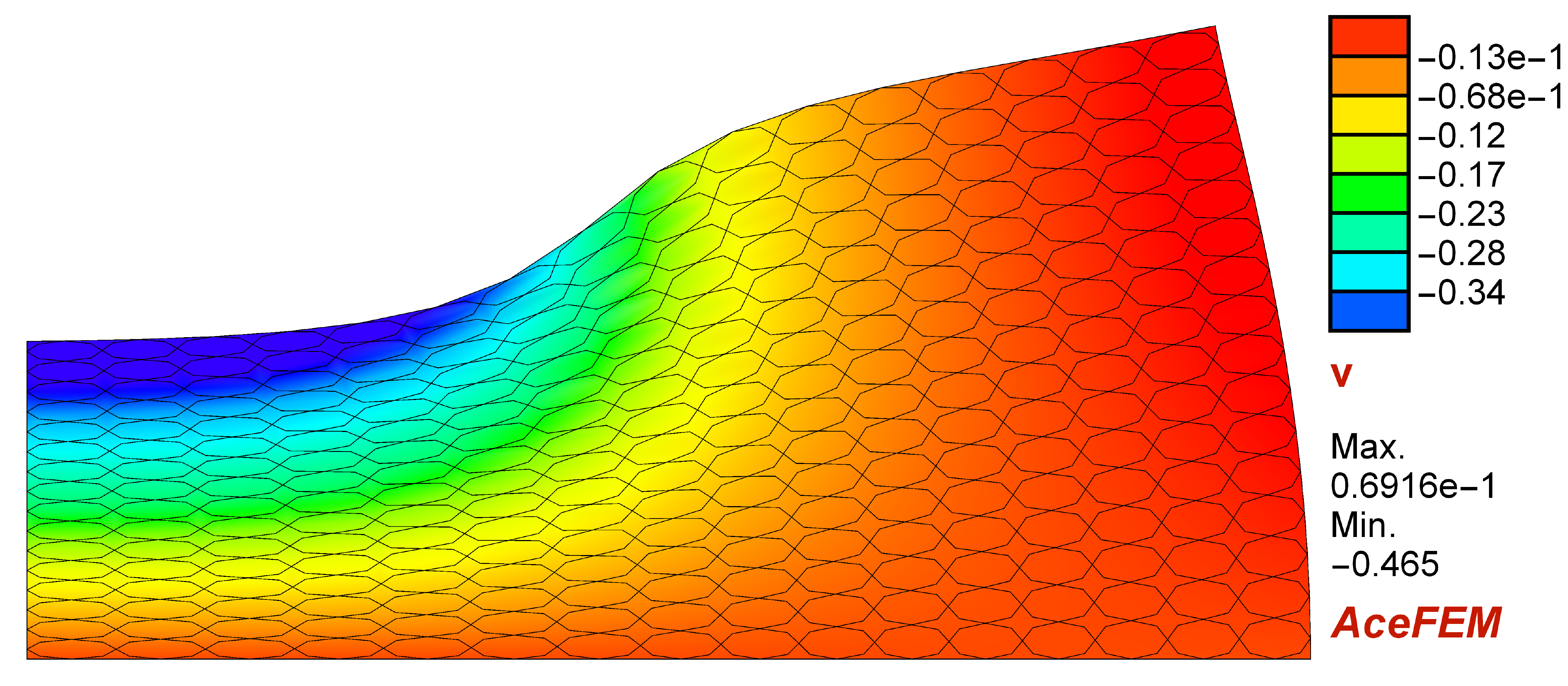}
		\caption{}
	\end{subfigure}
	\caption{Punch problem (a) geometry and (b) sample deformed configuration for an Ogden material model with $N=4$ and $\nu=0.3$ for an S\&S mesh. \label{fig:Punch}}
\end{figure}

\paragraph*{Punch problem} This problem consists of a rectangular body, vertically constrained along the bottom and horizontally constrained along the left and top faces, subjected to a uniformly distributed load $q_{\rm p}=200~\frac{\rm N}{\rm m}$ along half of the top edge, as depicted in Figure \ref{fig:Punch}(a). For this problem we consider an Ogden material model \cite{Ogden1997} with the strain energy function 
\begin{equation}
\Psi_{\rm Ogden} = \sum_{i=1}^{n} \frac{\mu_{i}}{\alpha_{i}} \left( \bar{\lambda}_{1}^{\alpha_{i}} + \bar{\lambda}_{2}^{\alpha_{i}} +  \bar{\lambda}_{3}^{\alpha_{i}} - 3 \right) + \frac{\kappa}{2}(J-1)^{2} \, , \label{eqn:OgdenComp}
\end{equation}
expressed in terms of the isochoric components of the principal stretches, see (\ref{eqn:iso1})-(\ref{eqn:iso3}), subject to the relation 
\begin{equation}
\mu = \sum_{i=1}^{n}\frac{\alpha_{i} \mu_{i}}{2} \, . \label{eqn:OgdenRatio}
\end{equation}
We consider a three term Ogden model and choose
\begin{eqnarray}
\alpha_{1}=1.3 \quad \text{and} \quad \alpha_{2}=5 \quad \text{and} \quad \alpha_{3}=-2 \, ,
\end{eqnarray}
and, for simplicity, define the ratios
\begin{eqnarray}
\mu_{1}=0.77 \mu  \quad \text{and} \quad \mu_{2}=0.1 \mu \quad \text{and} \quad \mu_{3}=-0.25 \mu  \, .
\end{eqnarray}
We choose a value of Young's modulus of $E_{y}=200~\rm Pa$. Figure \ref{fig:Punch}(b) shows the deformed configuration of the punch problem with an S\&S mesh for $N=4$ and a Poisson's ratio of $\nu =0.3$.
\\

\begin{figure}[ht!]
	\centering
	\begin{subfigure}{0.5\textwidth}
		\centering
		\def\svgwidth{\textwidth}
		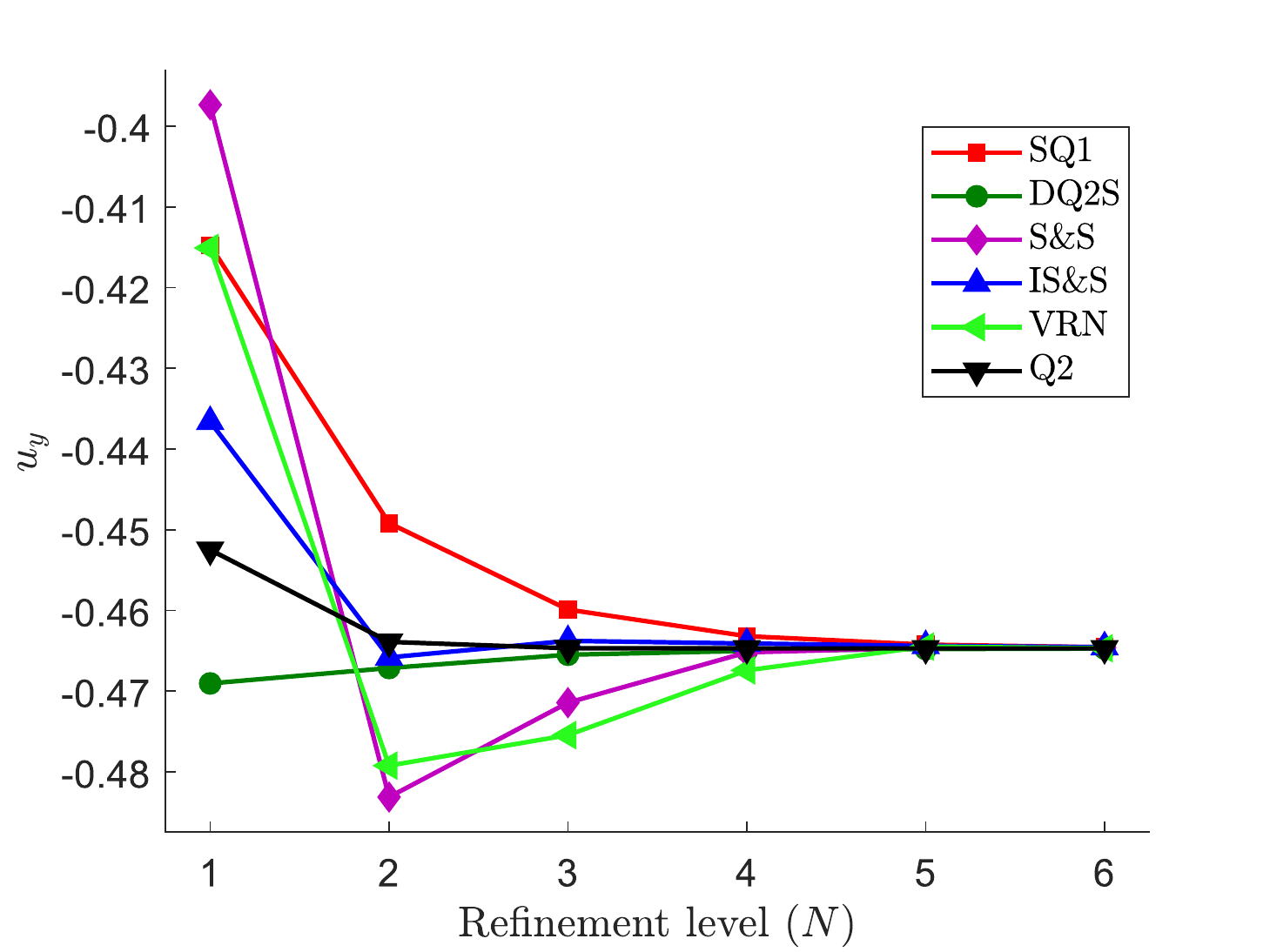
		\vspace*{-5mm}
		\caption{}
	\end{subfigure}%
	\begin{subfigure}{0.5\textwidth}
		\centering
		\includegraphics[width=\textwidth]{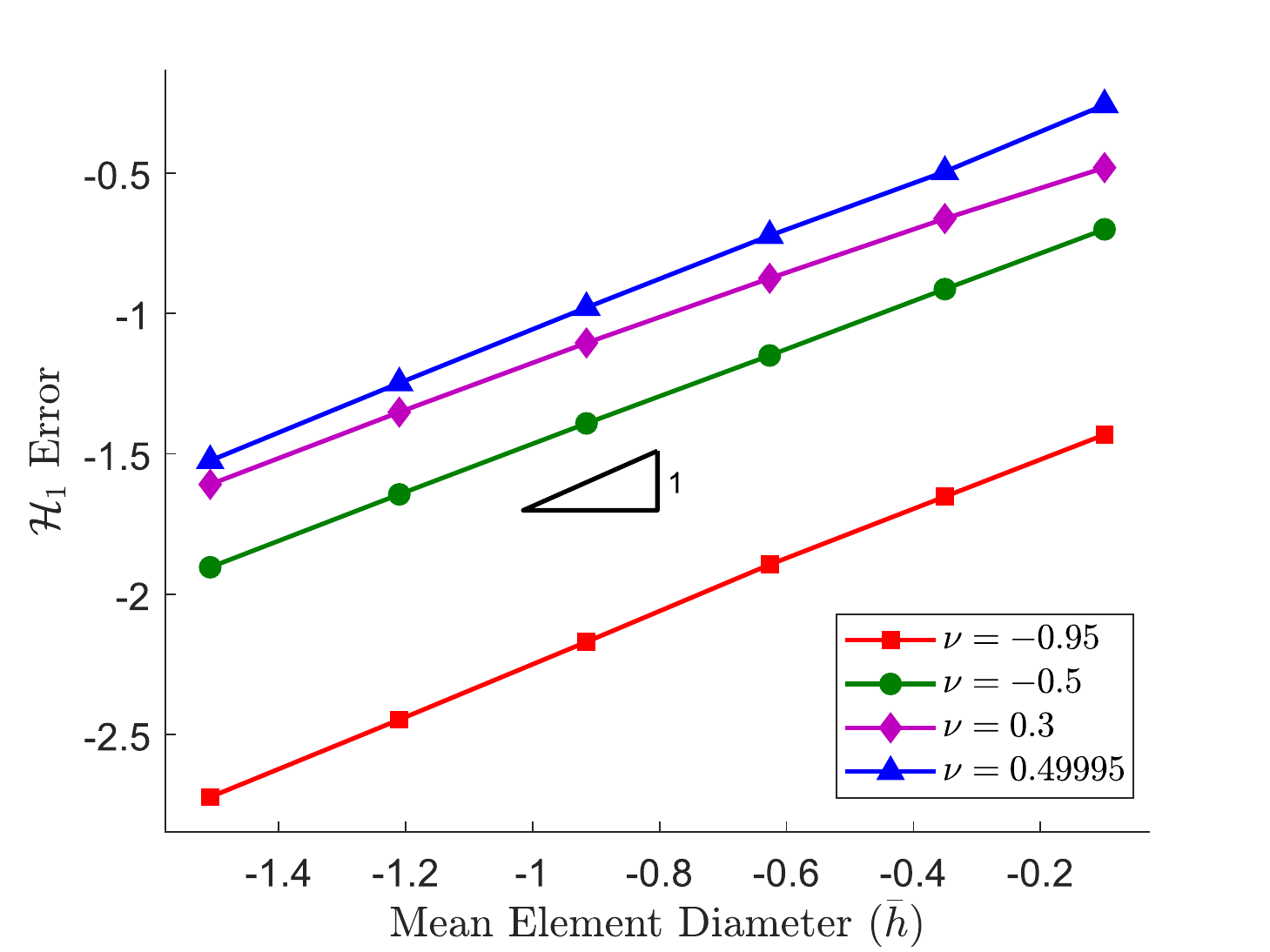}
		\vspace*{-5mm}
		\caption{}
	\end{subfigure}
	\caption{\normalsize Punch problem (a) tip deflection vs mesh refinement for $\nu=0.3$ and (b) loglog $\mathcal{H}_{1}$ error vs mean element diameter with S\&S meshes, for an Ogden material model. \label{fig:PunchConvergence}} \vspace*{-2mm}
\end{figure}

Figure \ref{fig:PunchConvergence}(a) shows the vertical displacement of the upper left hand corner of the body vs mesh refinement level with a Poisson's ratio of $\nu =0.3$. For S\&S and VRN meshes the displacement is initially erratic, due to the nature of the boundary conditions and the coarse meshes, but for $N\geq3$ the convergence of all formulations is smooth and monotonic, with good accuracy observed for meshes with $N\geq5$. Figure \ref{fig:PunchConvergence}(b) shows a plot of $\mathcal{H}_{1}$ error vs mean element diameter $\bar{h}$ on a loglog scale for the punch problem with S\&S meshes for a variety of choices of Poisson's ratio. We note for all choices of Poisson's ratio a near-optimal convergence rate with a gradient of approximately 1, with no suboptimal behaviour for coarse meshes with $\nu=-0.95$.

\begin{figure}[htp!]
	\centering
	\begin{subfigure}{.5\textwidth}
		\centering
		\includegraphics[width=\linewidth]{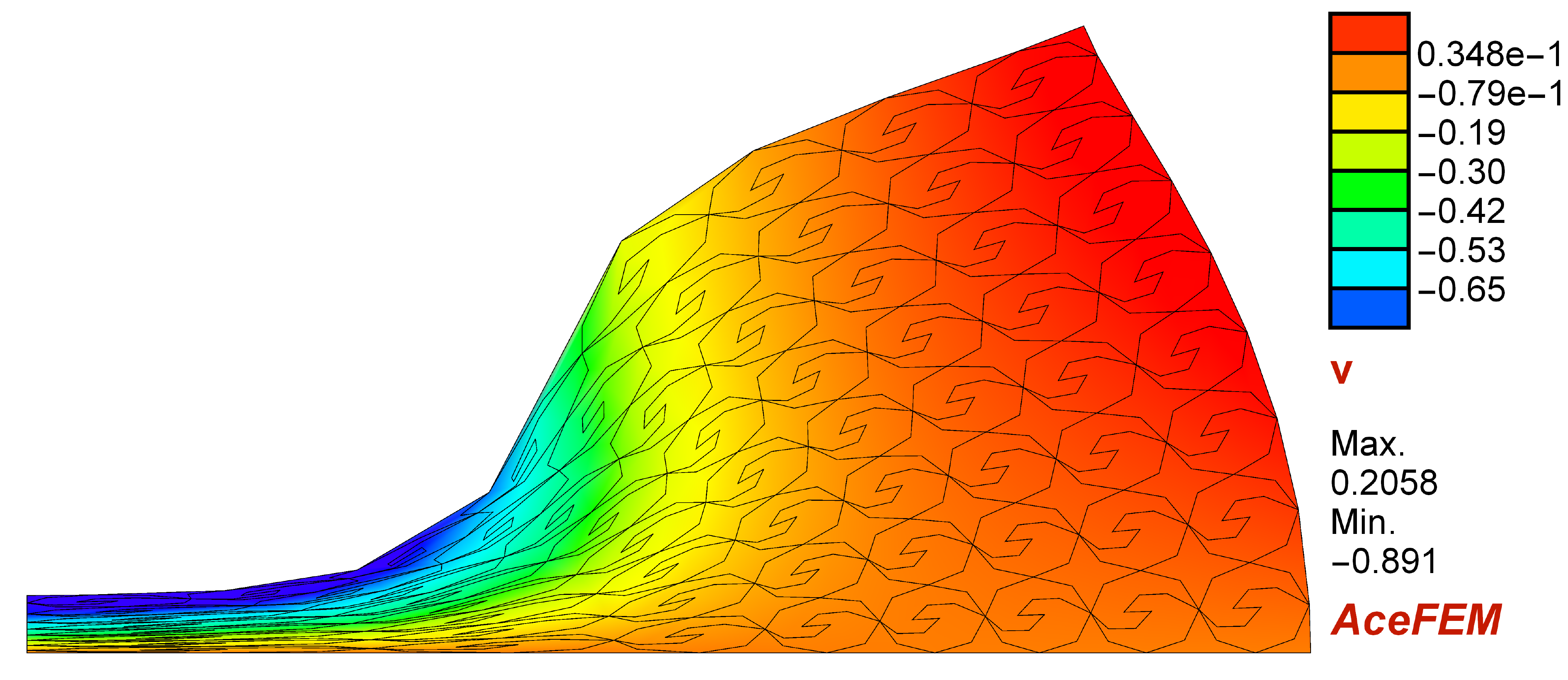}
		\vspace*{-5mm}
		\caption{}
	\end{subfigure}%
	\begin{subfigure}{.5\textwidth}
		\centering
		\includegraphics[width=\linewidth]{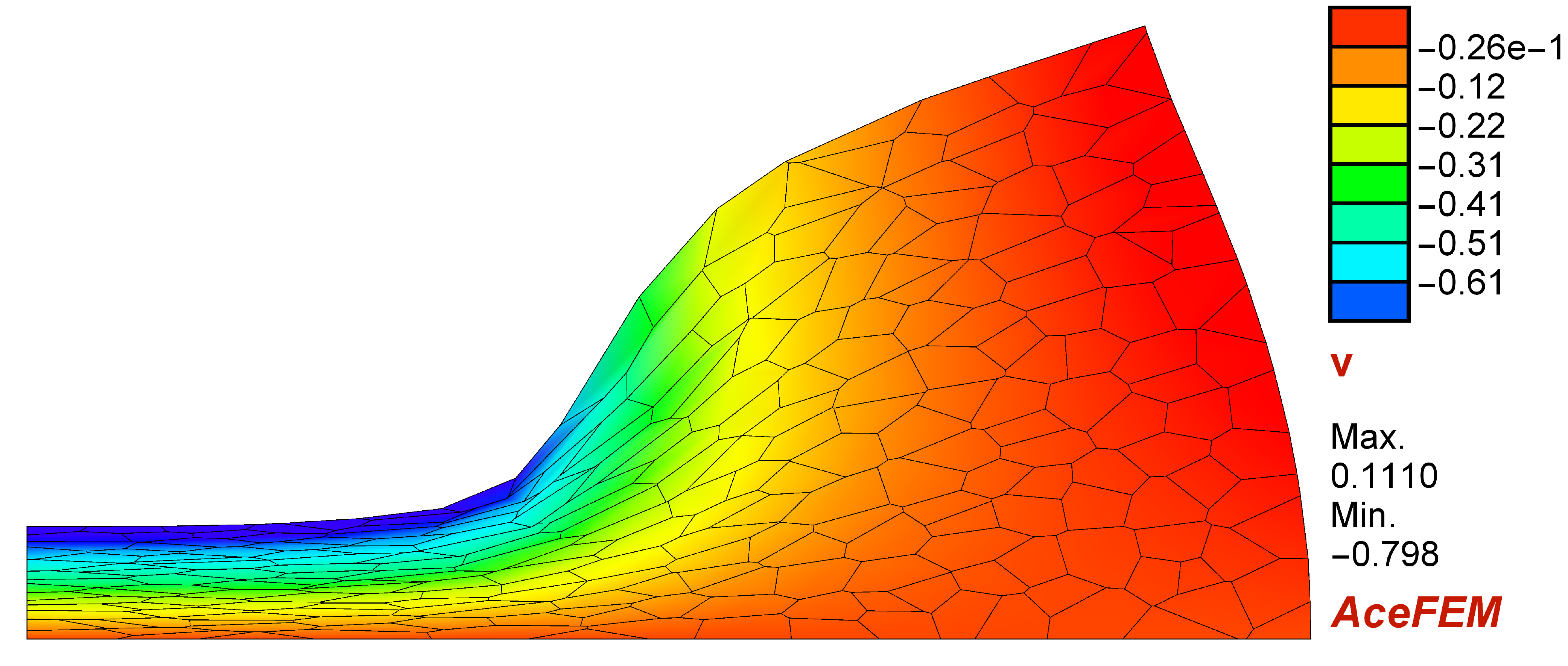}
		\vspace*{-5mm}
		\caption{}
	\end{subfigure}
	\caption{Deformed configuration of the punch problem with $q_{\rm p}=1000~\frac{\rm N}{\rm m}$ and $\nu=0.3$ for (a) a neo-Hookean material model with an IS\&S mesh and $N=3$, and (b) a Mooney-Rivlin material model with a VRN mesh and $N=4$. \label{fig:SeverePunch}}
\end{figure}

For all results thus far the loads have been chosen somewhat conservatively to ensure convergent solutions for all choices of Poisson's ratio. To demonstrate the ability of the VEM to undergo severe deformation we present in Figure \ref{fig:SeverePunch} the deformed configuration of the punch problem with $q_{\rm p}=1000~\frac{\rm N}{\rm m}$ and $\nu=0.3$. In Figure \ref{fig:SeverePunch}(a) we consider a neo-Hookean material model with an IS\&S mesh and $N=3$, while in \ref{fig:SeverePunch}(b) we consider a Mooney-Rivlin material model and a VRN mesh with $N=4$. The nature of the deformation is as expected, with the concave and interlocking elements showing robust behaviour under large deformations.

\begin{figure}[htp!]
	\centering
	\begin{subfigure}{.5\textwidth}
		\centering
		\includegraphics[width=\linewidth]{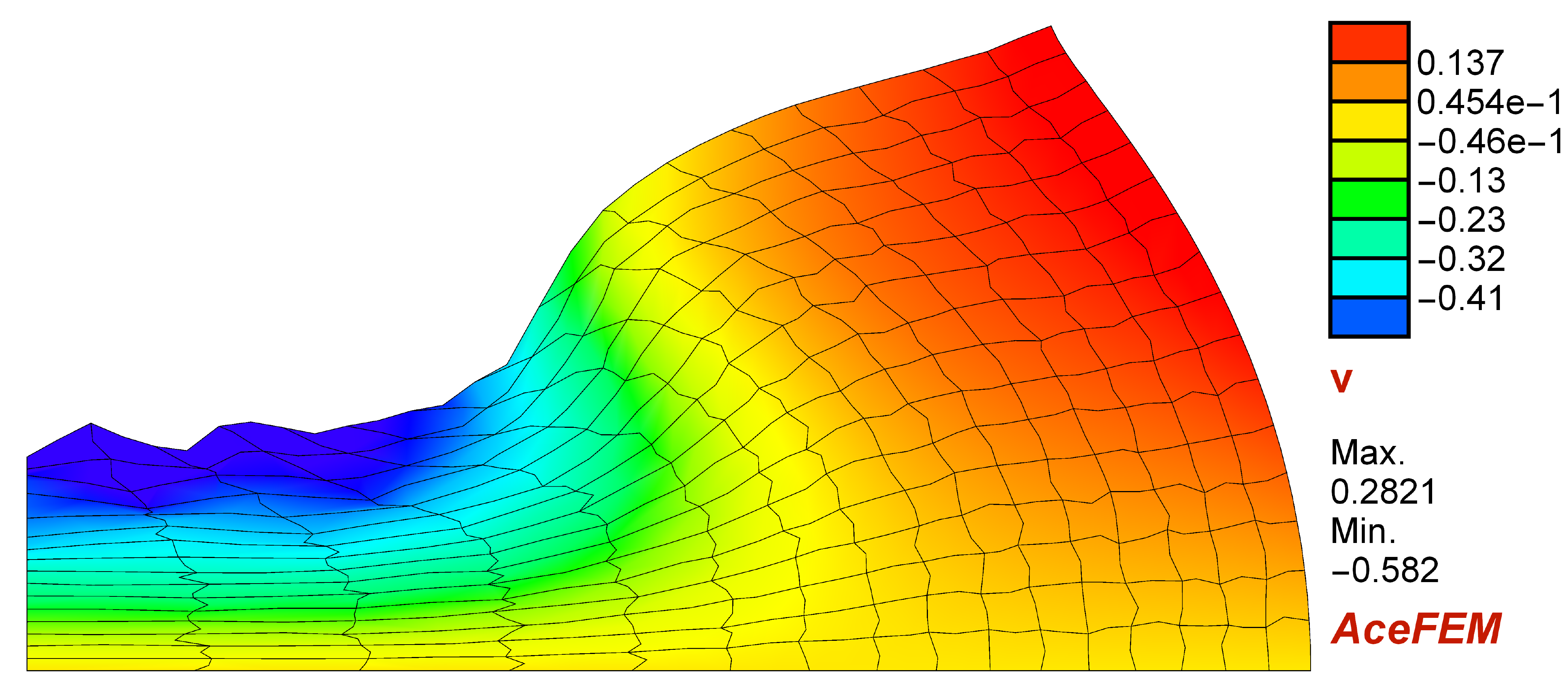}
		\vspace*{-5mm}
		\caption{}
	\end{subfigure}%
	\begin{subfigure}{.5\textwidth}
		\centering
		\includegraphics[width=\linewidth]{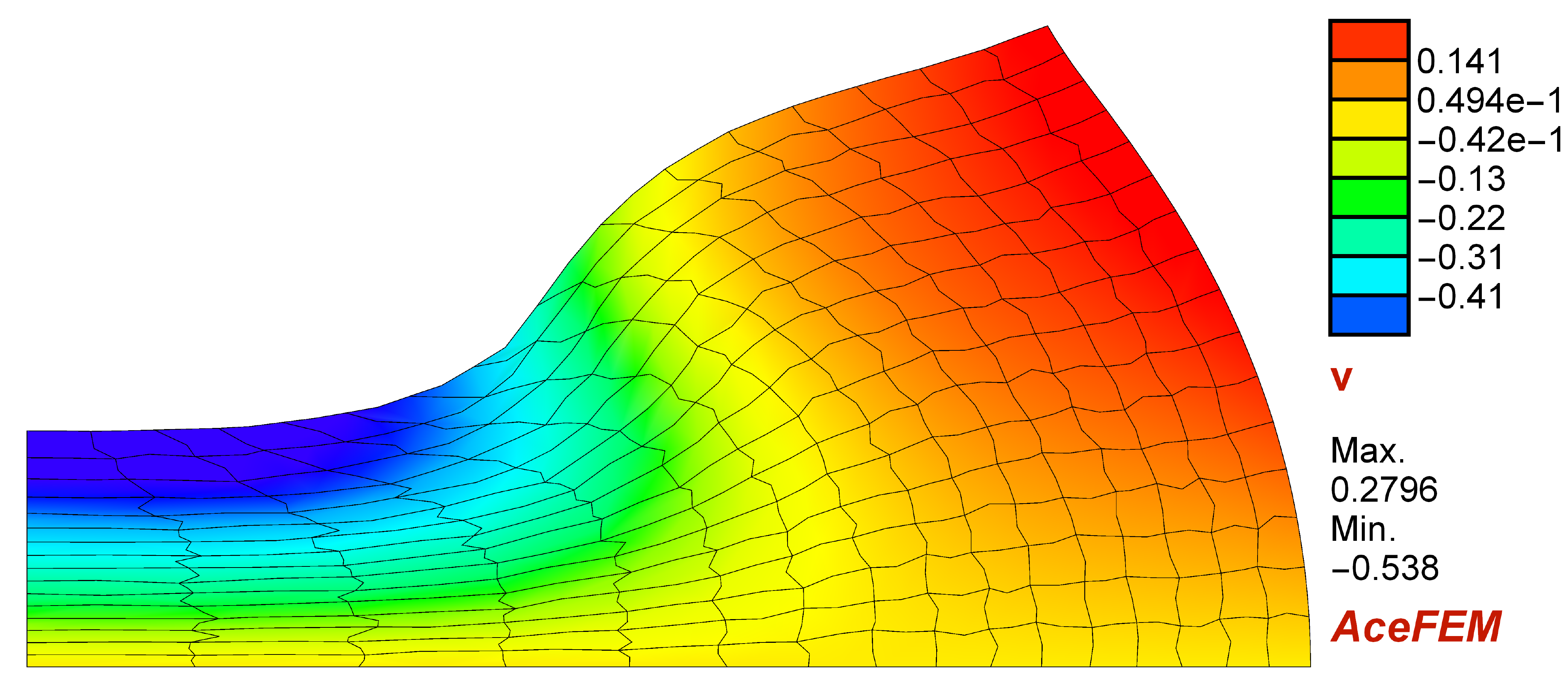}
		\vspace*{-5mm}
		\caption{}
	\end{subfigure}
	\caption{Deformed configuration of the punch problem for an Ogden material with $q_{\rm p}=1000~\frac{\rm N}{\rm m}$ and $\nu=0.49995$ for a DQ2S mesh with $N=4$ for (a) $\alpha=0$ and (b) $\alpha=\TEfive/E_{y}$. \label{fig:SeverePunchBoolOff}}
\end{figure}

In the results presented thus far the merit of the incompressibility term $\alpha$ has not been explicit. In Figure \ref{fig:SeverePunchBoolOff} we present the deformed configuration of the punch problem for the Ogden material model with $q_{\rm p}=1000~\frac{\rm N}{\rm m}$ and $\nu=0.49995$ for two cases. In Figure \ref{fig:SeverePunchBoolOff}(a) we set $\alpha=0$, `switching off' the incompressibility factor, and in Figure \ref{fig:SeverePunchBoolOff}(b) the $\alpha$-term is as described in Section \ref{subsec:Stab}. In Figure \ref{fig:SeverePunchBoolOff}(a) we note spurious, hourglass-like, deformation in the region of the applied load. The inclusion of the incompressibility factor, in Figure \ref{fig:SeverePunchBoolOff}(b), provides additional stabilization energy eliminating the spurious deformation and produces a more feasible deformation. Additionally, we again note that the VEM is locking free in the case of near-incompressibility.

\section{Concluding Remarks} \label{sec:Conclusion}
In this work we have formulated and implemented a virtual element method suitable for general use for plane isotropic hyperelasticity problems with a novel approach to the scaling of the Lam\'{e} parameters in the stabilization term. A variety of strain energy functions were investigated; for each of these, varying degrees of compressibility, including near-incompressibility, were considered under a range of different boundary conditions and loading types to test shear, bending and compression dominated deformations.

The virtual element method is found to be locking free in the case of near-incompressibility without any modification to the formulation. The locking-free behaviour is increasingly reported in the literature and is consistent with that observed in \cite{Chi2017,WriggersIsotropic2017,Reddy2019,Bellis2019}, to name a few. The convergence behaviour of the virtual element method in an $\mathcal{H}_{1}$-like norm is largely as expected, exhibiting near-optimal convergence rates. However, it was observed in the case of Cook's membrane problem that the virtual element method demonstrated unexpectedly high accuracy with a convergence rate higher than that observed for other numerical tests. Furthermore, with the inclusion of an incompressibility term in the stabilization, the virtual element method is found to be robust under severe deformations even with complex element geometries including concave as well as interlocking elements. 

Future work involves the extension of the presented stabilization methodology to anisotropic materials. It would be of interest to investigate the effects of lower order Taylor expansions of the $\hat{\lambda}$ term, as well as different choices of $\nu_{0}$ in the stabilization term. Of equal interest would be the extension of the work presented here to include higher order elements as well as problems in three dimensions. Additionally, formulation of a stabilization methodology not requiring small load steps or decomposition of a virtual element into triangles would represent a significant advancement to the virtual element method for non-linear problems.

\paragraph*{Acknowledgements} This work was carried out with support from the National Research Foundation of South Africa, through the South African Research Chair in Computational Mechanics. The authors acknowledge with thanks this support.




\bibliographystyle{elsarticle-num}
\bibliography{Isotropic_References}
\end{document}